\documentclass[12pt]{amsart}
\usepackage{graphicx}

\newtheorem{theorem}{Theorem}[section]

\newtheorem{lemma}[theorem]{Lemma}

\theoremstyle{definition}

\theoremstyle{remark}
\newtheorem{remark}[theorem]{Remark}

\numberwithin{equation}{section}

\begin{document}

\title[Dehn surgeries on $2$-bridge links] {
Dehn surgeries on $2$-bridge links\\
which yield reducible $3$-manifolds}

\author{Hiroshi Goda, Chuichiro Hayashi and Hyun-Jong Song}


\thanks{
This work was supported by Joint Research Project
$\lq$Geometric and Algebraic Aspects of Knot Theory',
under the Japan-Korea Basic Scientific Cooperation Program
by KOSEF and JSPS. The authors would like to thank
Professor Hitoshi Murakami for giving us this opportunity.
The first and second authors are partially supported
by Grant-in-Aid for Scientific Research,
(No. 15740031 and No. 15740047 respectively),
Ministry of Education, Science, Sports and Technology, Japan.
}

\maketitle


\section{Introduction}

 No surgery on a non-torus $2$-bridge knot yields a reducible $3$-manifold
as shown in Theorem 2(a) in \cite{HT} by A. Hatcher and W. Thurston.
 Dehn surgeries on $2$-bridge knots are already well studied
by M. Brittenham and Y.-Q. Wu in \cite{BW}.
 See also \cite{P}.

 On $2$-bridge links of $2$-components,
 Wu showed in \cite[Theorem 5.1 and Remark 5.5]{W} the following theorem.
 The universal covering space of a laminar $3$-manifold $M$
is the Euclidean $3$-space ${\Bbb R}^3$,
and then $M$ is not reducible by \cite{GO}.

\begin{theorem}\label{thm:Wu}
{\rm (\cite{W})}
 A non-trivial Dehn surgery on a $2$-bridge link yields a laminar $3$-manifold
except for the $2$-bridge links $L([r,s])$ with $[r,s]=1/(r+1/s)$.
 For them with $r,s \ne \pm 1$,
we obtain laminar $3$-manifolds
if the two surgery slopes are both non-integral.
\end{theorem}

 In this paper we study Dehn surgery on $2$-bridge links $L([r,s])$.
 Note that $L([r,s])$ is a link of $2$-components
if and only if both $r$ and $s$ are odd integers.
 Since $L([r,s]) \cong L([-s,-r])$
and $L([-r,-s]) \cong L([s,r])$ is the mirror image of $L([r,s])$,
it is enough to consider $L([r,s])$ with $r > 0$.
 When $|r|=1$ or $|s|=1$,
$L([r,s])$ is a $(2, k)$-torus link for some even integer $k$.
 Thus we consider the case where $r \ge 3$ and either $s \ge 3$ or $s \le -3$.

\begin{figure}[htbp]
\begin{center}
\includegraphics[width=3cm]{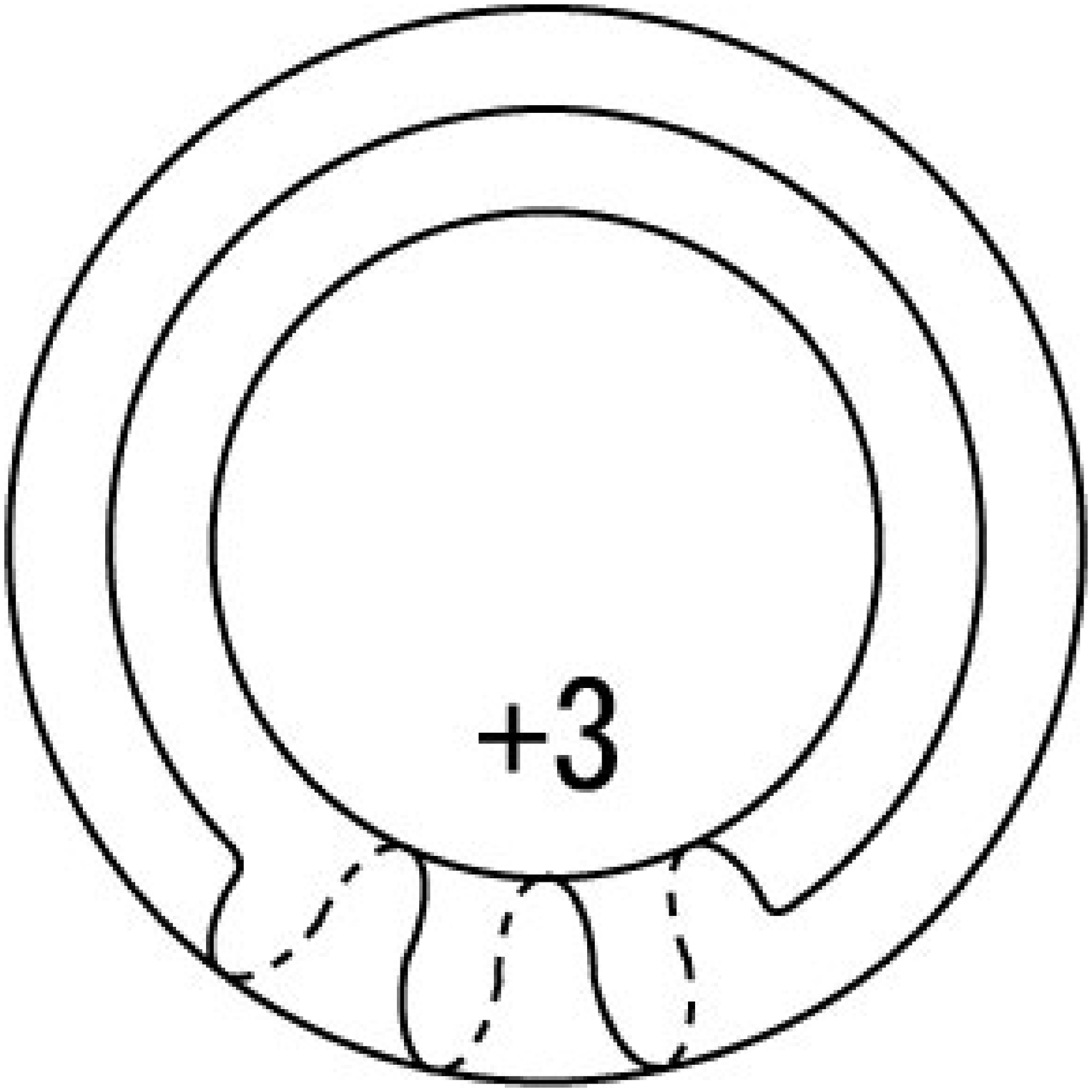}
\end{center}
\caption{}
\label{fig:SurgSlop}
\end{figure}

 Every component of a $2$-bridge link is the trivial knot.
 We coordinate slopes on the toral boundaries of the link exterior
so that the slope in Figure \ref{fig:SurgSlop} is $+3$ rather than $-3$.
 Let $L(p/q)[\gamma_1, \gamma_2]$
be the $3$-manifold obtained by the Dehn surgery on the $2$-bridge link $L(p/q)$
with surgery slopes $\gamma_1$ and $\gamma_2$.
 The order of the surgery slopes $\gamma_1, \gamma_2$ is not crucial
since every $2$-bridge links has a symmetry of $\pi$-rotation
which exchanges the two components.
 See Figure \ref{fig:3over8}.

\begin{theorem}\label{thm:red}
 $L([2w+1,2u+1])[\gamma_1, \gamma_2]$
with $w \ge 1$ and either $u \ge 1$ or $u \le -2$
is a reducible $3$-manifold
if and only if it is one of the followings or its mirror image.
\item\begin{enumerate}
\item
\begin{enumerate}
\item
 $\gamma_1 = -w+u+1$ and $\gamma_2 = -w+u-1$.
\item
 $\gamma_1 = -w+u$ and $\gamma_2 = -w+u$.
\end{enumerate}
\item
 $w=1$, $u=-2$ and $\gamma_1 = -1$ and $\gamma_2 = -6$.
\end{enumerate}
\end{theorem}

 If we perform Dehn surgery of slope $\gamma$
on one component of a $2$-bridge link $L(p/q)$,
then the other component forms a knot
in a lens space, $S^1 \times S^2$ or $S^3$.
 We denote it by $L(p/q)[\gamma]$.
 The next two theorems show
when it is a torus knot or a cable knot in $S^3$.
 More general version on knots in lens spaces is given
in Theorem \ref{thm:TorusCable}.

\begin{theorem}\label{thm:torus}
 Let $L(p/q)$ be a non-torus $2$-bridge link.
 Then $L(p/q)[\gamma]$ is a non-trivial torus knot in $S^3$
if and only if it is one of the followings or its mirror image.
\begin{enumerate}
\item
 $p/q=[p+2,p]$, $\gamma=1$
 and $L(p/q)[\gamma]$ is the $(p,-(p+2))$-torus knot.
\item
 $p/q=[3,3]$, $\gamma = 1$
 and $L(p/q)[\gamma]$ is the $(2,-5)$-torus knot.
\item
 $p/q=[-3,3]$, $\gamma = - 1$
 and $L(p/q)[\gamma]$ is the $(2,-3)$-torus knot.
\end{enumerate}
\end{theorem}

\begin{theorem}\label{thm:cable}
 Let $L(p/q)$ be a non-torus $2$-bridge link.
 Then $L(p/q)[\gamma]$ is a cable knot in $S^3$
if and only if it is one of the followings or its mirror image.
\begin{enumerate}
\item
 When $w=u$, $L([2w+1, 2u+1])[1]$
is the $(2, -2w^2-2w-1)$-cable of the $(w, -w-1)$-torus knot.
\item
 When $w=u+2$, $L([2w+1, 2u+1])[-1]$
is the $(2, 2w^2-2w-1)$-cable of the $(-w, -w+1)$-torus knot.
\end{enumerate}
\end{theorem}


\begin{remark}
\noindent
\begin{enumerate}
\item
 Theorem \ref{thm:torus} (1) gives
counter examples for Ait-Nouh and Yasuhara's conjecture:
if a $(p,q)$-torus knot ($q \ge p > 0$) is obtained
by a twisting operation on the trivial knot,
then $q = np \pm 1$ for some integer $n$.
 See \cite{AY}.
 The existence of the counter examples has already been shown
in the previous paper \cite{GHS}.
\item
 The essential planar surface corresponding to Theorem \ref{thm:red} (3)
has $4$ boundary circles of slope $-1$ on a component of the link,
and $2$ boundary circles of slope $-6$ on the other component.
 See Calculation for $d_{26}$ in Section 10.
 This shows that the cabling conjecture for  hyperbolic knots in $S^3$
can not be extended to those in a lens space.
 Because $L([3,-3])[-6]$ is a hyperbolic knot in the $(6,-1)$-lens space
according to the computer software SnapPea programmed by J. Weeks. 
\item
 $L([3,3])[\pm 1]$ is the $(2, \mp 5)$-torus knot
as shown in Theorem \ref{thm:torus} (2).
 Note that $L([3,3])$ is amphicheiral.
\end{enumerate}
\end{remark}

 In general, let $M$ be a $3$-manifold,
and $K$ a knot in the interior of $M$.
 Then $K$ is called a {\it composite knot}
if there is a $2$-sphere $S$ intersecting $K$ transversely in two points in $M$
such that $S$ does not bound a $3$-ball which $K$ intersects in a trivial arc.
 Otherwise, $K$ is called a {\it prime knot}.
 We say that $K$ is a {\it satellite knot}
if its exterior $E$ contains an incompressible torus
which is not parallel to a toral boundary component of $E$.

\begin{theorem}\label{thm:satellite}
 Suppose $L(p/q)$ is a non-torus $2$-bridge link.
 If $L(p/q)[\gamma]$ is a prime satellite knot,
then either (1) $L(p/q) \cong L([2w+1, 2u+1])$
for some integer $w, u$ with $w \ge 2$ and either $u \ge 2$ or $u \le -3$,
and $\gamma = -w+u \pm 1$,
or (2) $L(p/q) \cong L([2w, v, 2u])$
for some integers $w, v, u$ with $|w|, |v|, |u| \ge 2$.
 Conversely, in case (1), $L([2w+1, 2u+1])(-w+u \pm 1)$ is a satellite knot.
\end{theorem}

 In fact, in case (1), it is a cable knot of a torus knot
in the $(-w+u \pm 1, 1)$-lens space
by Theorem \ref{thm:TorusCable}.

\vspace{2mm}

 It is well-known that,
if $L(p/q)[\gamma_1, \gamma_2]$ is reducible,
then the exterior of $L(p/q)$ contains an essential planar surface
with boundary slopes $\gamma_i$ on $\partial N(L_i)$
for $i=1$ and $2$,
where $\lq\lq$essential" means incompressible and boundary incompressible.

 Our study is based on the classification of essential surfaces
in $2$-bridge link exteriors by W. Floyd and A. Hatcher (\cite{FH}).
 We assume the readers' good familiarity with their study.
 The space of isotopy classes of the surfaces is much more complicated
than that of $2$-bridge knot (\cite{HT}).
 A. Lash gives a way to calculation of boundary slopes for such surfaces,
and did calculate for the Whitehead link $L([3,-3])$ in \cite{L}.
 We should note
that J.Hoste and P.Shanahan \cite{HS} improved the way of Lash recently.
 However, the way of calculation of genera of the surfaces is not given there.
 Note that the surfaces carried by branched surfaces given in \cite{FH}
may be non-orientable or disconnected.

 This paper is organized as follows.
 In Section 2,
we recall results in \cite{FH}, \cite{L}.
 In Section 3, we obtain all the minimal edge-paths for $L([r,s])$.
 In Section 5,
we give a way to calculate Euler characteristics of such surfaces.
 In Section 6,
we generalize the notion of genus for disconnected or non-orientable surfaces.
 The generalized genus is zero
if and only if there is either a planar component
or a projective plane with holes component.
 In Section 7 and 8,
we calculate boundary slopes, Euler characteristics, generalized genus
of the surfaces carried by the branched surfaces
corresponding to certain two minimal paths.
 In Section 9,
we give data of all the essential surfaces
in the exterior of $L([r,s])$.
 In Section 10 and Appendix A,
the proof of $\lq\lq$only if part" of Theorem 1.2 is given.
 In Section 11, we give the proofs of sufficiency of 
Theorems \ref{thm:red}, \ref{thm:torus}, \ref{thm:cable}
and \ref{thm:TorusCable}.
 In Section 12, we prove Theorem 1.6.

\section{A summary of results of Floyd, Hatcher and Lash}

 W. Floyd and A. E. Hatcher studied the spaces of incompressible surfaces
in $2$-bridge link exteriors in \cite{FH}.
 A. E. Lash studied the way how to compute
the spaces of boundary slopes of incompressible surfaces
for $2$-bridge links in Chapter 1 in his doctoral dissertation \cite{L}.
 (He calculated the space of boundary slopes
for the Whitehead link $L([3,-3])$.)

 In this section, we briefly recall their results.
 Very roughly speaking,
every orientable essential surface is carried
by a branched surface corresponding to a minimal edge-path
in a certain planar graph in the upper half plane,
as below.

{\bf The diagrams of slope system and minimal edge-paths}.

 The diagram $D_1$ is an embedded graph
on the upper half plane ${\Bbb H}$
with the real line ${\Bbb R}$ and the point at infinity $1/0$
such that its vertices are the rational points in ${\Bbb R} \cup \{ 1/0 \}$
and that its edges are geodesics on the upper half model of the hyperbolic plane
and join two vertices $a/c$, $b/d$ if and only if $ad-bc = \pm 1$.
 We regard the vertex $0$ as $0/1$.
 See Figure 1.1 in \cite{FH},
where the diagram $D_1$ transformed onto the Poincar\'{e} disc model
by $-\dfrac{z-\frac{1+i}{2}}{z-\frac{1-i}{2}}$ is described.
 (Figure 4 in \cite{HT} is transformed by $-\dfrac{z-i}{z+i}$.)

 We identify the matrix
$\left( \begin{array}{cc}
a & b \\
c & d \\
\end{array} \right)
\in PSL_2 {\Bbb R}$
with the M\"{o}bius transformation $\dfrac{az+b}{cz+d}$,
an orientation preserving isometry
of the upper half model of the hyperbolic plane.
 Let $G \subset PSL_2({\Bbb Z})$ be the subgroup of transformations
$\left( \begin{array}{cc}
a & b \\
c & d \\
\end{array} \right)$
with $c$ being an even integer.
 $G$-images of the ideal quadrilateral $\langle 0/1, 1/0, 1/1, 1/2 \rangle$ form
a tiling of the upper half plane by quadrilaterals.
 These edges form the $G$-orbit of the edge $\langle 1/0, 0/1 \rangle$
and we label them by $A$.
 The other edges of $D_1$ are diagonals of the quadrilaterals
and form the $G$-orbit of the edge $\langle 0/1, 1/1 \rangle$.
 We label them by $C$.

 We form the diagram $D_0 = D_{\infty}$ from $D_1$
by deleting all the edges labeled $C$
and adding the $G$-orbit of the edge $\langle 1/0, 1/2 \rangle$.
 The new edges are the opposite diagonals of the quadrilaterals,
and we label them by $D$.
 The edges of the quadrilaterals are labeled by $B$
in the diagram $D_0 = D_{\infty}$.

 For $0 < t < \infty$ with $t \ne 1$,
the diagram $D_t$ is obtained from $D_1$
by deleting the diagonal edge labeled $C$ in each quadrilateral $Q$
and adding a small rectangle having a vertex in the interior of each edge of $Q$
so that $g(D_t) = D_t$ for arbitral $g \in G$.
 (We distinguish these small rectangles from the quadrilaterals
by keeping to call them rectangles rather than quadrilaterals,
to avoid confusion.)
 The edge $\langle 1/0, 0/1 \rangle$ is divided to two edges by the added vertex.
 We label by $A$ the $G$-orbit of one of them including the vertex $1/0$,
and that of the other including $0/1$ by $B$.
 We label by $C$ (resp. $D$) an edge of an added small rectangle
if it cobounds a triangle face
together with two edges with label $A$ (resp. $B$).
 The edges of $D_t$ is fall into four $G$-orbits
corresponding to these four labels.

 For a given reduced rational number $p/q$,
an oriented edge-path $\lambda$ from $1/0$ to $p/q$ in $D_t$
with $0 \le t \le \infty$
is called {\it minimal}
if no two consecutive edges of $\lambda$ lie
in the boundary of the same triangle face or rectangle face of $D_t$.
 There is a unique finite sequence of quadrilaterals $Q_{p/q}$
such that the first one contains the vertex $1/0$,
that the last one contains the vertex $p/q$,
and that every pair of consecutive ones intersect in a single edge.
 Every minimal edge-path from $1/0$ to $p/q$
is entirely contained in $Q_{p/q}$.
 Hence there are only finitely many such paths for $p/q$.

\begin{figure}[htbp]
\begin{center}
\includegraphics[width=6cm]{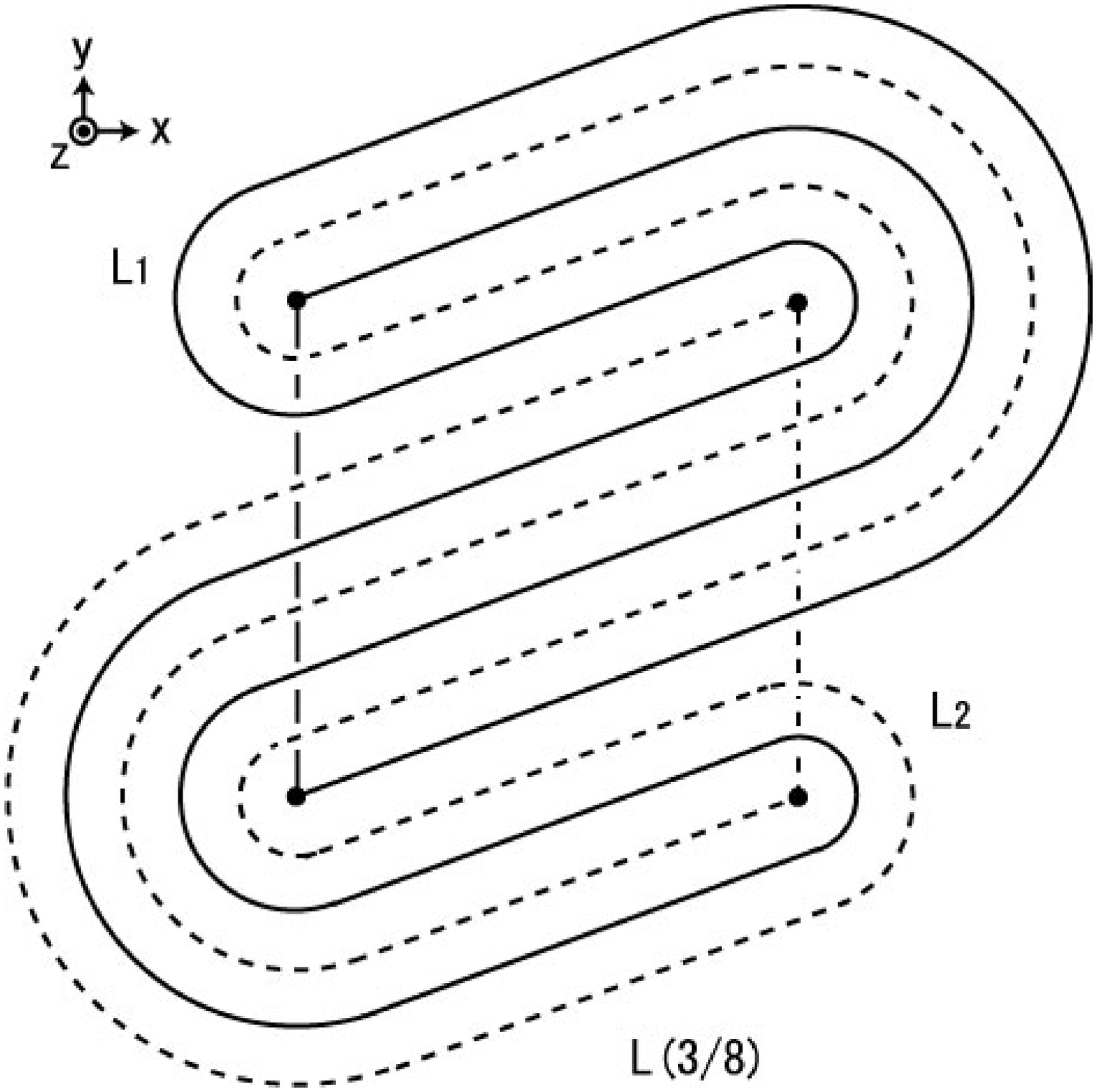}
\end{center}
\caption{}
\label{fig:3over8}
\end{figure}

{\bf Links and branched surfaces}.

 The $3$-sphere $S^3$ can be regarded
as the $2$-points compactification of $S^2 \times {\Bbb R}$,
and each level sphere $S^2 \times \{ r \}$
as the quotient ${\Bbb R}^2 / \Gamma$,
where $\Gamma$ is the group
generated by $180^{\circ}$ rotations about integer points in ${\Bbb R}^2$.
 The $\Gamma$ image of a subset $S \subset {\Bbb R}^2$
is denoted by $\Gamma(S)$.
 Then such an $S^2$ has precisely four integer points
$\Gamma(0,0)$, $\Gamma(1,0)$, $\Gamma(0,1)$ and $\Gamma(1,1)$.
 For a reduced fraction $p/q$ with $p, q \in {\Bbb Z}$,
the $2$-bridge link $L(p/q)$ can be placed in $S^2 \times [0,1]$
so that
(1) $L(p/q) \cap (S^2 \times \{ 0 \})$
is the union of the two arcs of slope $\infty$
with their endpoints in the integer points,
(2) $L(p/q) \cap (S^2 \times \{ r \})$, $r \in (0,1)$
consists of the four integer points and
(3) $L(p/q) \cap (S^2 \times \{ 1 \})$
is the union of the two arcs of slope $p/q$
with their endpoints in the integer points.
 In $S^2 = {\Bbb R}^2 / \Gamma$ with the four integral points,
the slope of the arcs are determined
by the slopes of the preimage lines in ${\Bbb R}^2$.
 $L(p/q)$ is a $2$-component link rather than a knot
if and only if $q$ is an even integer.
 This positioning shows that $L(p/q)$
has a symmetry by the $180^{\circ}$ rotation
about the axis $\Gamma(1/2,1/2) \times [0,1]$,
which exchanges the two components of $L(p/q)$.
 See Figure \ref{fig:3over8}.

 $SL_2 ({\Bbb Z})$ acts on ${\Bbb R}^2$ by
$\left( \begin{array}{cc}
a & b \\
c & d \\
\end{array} \right)
\left( \begin{array}{c}
y \\
x \\
\end{array} \right)
=
\left( \begin{array}{c}
ay + bx \\
cy + dx \\
\end{array} \right)$.
 (This act is in a little unusual manner
because we write the first coordinate under the second
in the above two $2$-dimensional vectors.)
 Hence $PSL_2 ({\Bbb Z})$ acts
on the level sphere $S^2 \times \{ r \} = {\Bbb R}^2 / \Gamma$.
 $\left( \begin{array}{cc}
a & b \\
c & d \\
\end{array} \right)$
brings the slope $p/q$ to the slope $\dfrac{ap+bq}{cp+dq}$,
corresponding to the M\"{o}bius transformation $\dfrac{az+b}{cz+d}$
on the upper half plane.

 The vertices of the diagrams $D_1, D_0 = D_{\infty}, D_t$
correspond to the slopes of arcs in the level spheres.
 For every minimal edge-path $\gamma$ in $D_t$ with $1 < t < \infty$,
we construct a corresponding branched surface $\Sigma_{\gamma}$ as below.
 Let $e_1, \ldots, e_k$ be the sequence of edges of $\gamma$.
 For each edge $e_i$, the sub-branched surface
$\Sigma_{e_i} = \Sigma_{\gamma} \cap (S^2 \times [(i-1)/k, i/k])$
is an adequate transformation
of the branched surface $\Sigma_A, \Sigma_B, \Sigma_C, \Sigma_D$
in Figure 3.1 in \cite{FH}
according to the label of $e_i$.
 (In the figure, 
we can find the four strand
$\Gamma(0,0), \Gamma(1,0), \Gamma(0,1), \Gamma(1,1)$ of a link
from the left to the right.)
 Precisely,
$\Sigma_A, \Sigma_B, \Sigma_C, \Sigma_D$ is the branched surface
corresponding to the edges $A_0, B_0, C_0, D_0$ as below.
 $A_0$ is oriented from $1/0$ to the interior point of $\langle 1/0, 0/1 \rangle$,
$B_0$ from $0/1$ to the interior point of $\langle 1/0, 0/1 \rangle$,
$C_0$ from the interior point of $\langle 1/0, 1/1 \rangle$
to the interior point of $\langle 1/0, 0/1 \rangle$ and
$D_0$ from the interior point of $\langle 0/1, 1/2 \rangle$
to the interior point of $\langle 1/0, 0/1 \rangle$.
 (In fact, the arcs at the bottom of $\Sigma_A$ are of slope $1/0$
and the arcs at the top of $\Sigma_A$ are of slopes $0/1$ and $1/0$.)
 There is a M\"{o}bius transformation $g \in G$
with $e_i = g(e_0)$, $e_0 \in \{ A_0, B_0, C_0, D_0 \}$
according to the label of $e_i$.
 If the orientations of $e_i$ and $g(e_0)$ match,
then $\Sigma_{e_i} = (g \times *) \Sigma_{e_0}$,
where $*$ is the map $[0,1] \ni s \mapsto (i+s-1)/k \in [(i-1)/k, i/k]$.
 If the orientations don't match, we reflect $\Sigma_{e_i}$ upside down.
 Rotating $\Sigma_{\gamma}$
by $180^{\circ}$ about $\Gamma(1/2,1/2) \times [0,1]$,
we obtain the branched surface for $D_{1/t}$,
where $0< 1/t < 1$.

 The diagram $D_{\infty}$
has only edges labeled $B$ and edges labeled $D$.
 We obtain $\Sigma_B, \Sigma_D$ for $D_{\infty}$
from $\Sigma_B, \Sigma_D$ in Figure 3.1 in \cite{FH}
by deleting square sectors (arcs of slope $0) \times [0,1]$ labeled $\beta$.
 The resulting $\Sigma_B$ (resp. $\Sigma_D$) corresponds
to the edge from $0/1$ to $1/0$ (resp. from $1/2$ to $1/0$).
 Rotating $\Sigma_B, \Sigma_D$ for $D_{\infty}$
by $180^{\circ}$ about $\Gamma(1/2,1/2) \times [0,1]$,
we obtain those for $D_0$.

 The edges of $D_1$ are labeled $A$ or $C$.
 We obtain $\Sigma_A$ for $D_1$ from $\Sigma_A$ in Figure 3.1 in \cite{FH}
by deleting the square (the arc of slope $1/0) \times [1/2,1]$
labeled $\alpha - \beta$.
 The resulting $\Sigma_A$ corresponds to the edge from $1/0$ to $0/1$.
 $\Sigma_C$ for $D_1$ is obtained from $\Sigma_A$ for $D_1$
by transformed by
$\left( \begin{array}{cc}
1 & 0 \\
1 & 1 \\
\end{array} \right)$.
 This corresponds to the edge from $1/1$ to $0/1$.
 ($\Sigma_C$ for $D_1$ is obtained
also from $\Sigma_C$ in Figure 3.1 in \cite{FH}
by deleting the square (the arc of slope $1/0) \times [0,1]$
labeled $\alpha-\beta$
and adding a complementary horizontal square
of the horizontal square in the level sphere $S^2 \times \{ 1/2 \}$.)

 Then the branched surfaces for $D_0 = D_{\infty}$ and $D_1$ are constructed
similarly as those for $D_t$ with $0<t<\infty$, $t \ne 1$.

 Let $L_1$ (resp. $L_2$) be the component of $L(p/q)$
containing the arcs $\Gamma(0,0) \times [0,1]$ called strand 1
and strand 2 $\Gamma(0,1) \times [0,1]$
(resp. strand 3 $\Gamma(1,0) \times [0,1]$
and strand 4 $\Gamma(1,1) \times [0,1]$).
 For a surface $F$ in the exterior of $L(p/q)$,
let $\alpha$ (resp. $\beta$) be the minimal number of intersection points
of the boundary circles
$\partial F \cap \partial N(L_1)$ (resp. $\partial F \cap \partial N(L_2)$)
and a meridian circle of $\partial N(L_1)$.
 Then $t = \alpha / \beta$,
which is the subscript of $D_t$ with $0 \le t \le \infty$.
 Suppose that $F$ is carried by a branched surface $\Sigma_{\gamma}$
as above.
 In Figure 3.1 in \cite{FH},
where $1 < t < \infty$ and $\alpha > \beta$,
the labels $\beta$, $\alpha-\beta$, $(\alpha-\beta)/2$,
$n$, $\beta-n$ and $\alpha -\beta -n$
indicates the number of sheets carried by the sectors.
 The branching number $n$ is determined for each segment $\Sigma_{e_i}$ by $F$.
 $0 \le n \le \beta$ for $\Sigma_A$
and $0 \le n \le \alpha - \beta$ for $\Sigma_D$
when $1 \le t \le \infty$.

 As shown in Theorem 3.1 (a) in \cite{FH},
every $\lq\lq$essential" surface in the exterior of $L(p/q)$
is carried by some branched surface
corresponding to a minimal edge-path from $1/0$ to $p/q$
in $D_t$ with $t = \alpha / \beta$.
 Conversely, an orientable surface carried by such a branched surface
is essential.
 There, an {\it essential} surface is
an incompressible and meridionally incompressible surface
without peripheral component.
 (We don't need to consider meridionally incompressibility
because a surface with non-meridional boundary slope
both on $\partial N(L_1)$ and $\partial N(L_2)$
is always meridionally incompressible.)

 However, such a branched surface may carry
non-orientable or disconnected surfaces.
 (Moreover, there may be an essential non-orientable surface
which is not carried by any branched surface as above.)

{\bf Boundary slopes}.

 A. E. Lash calculated the space of boundary slopes for the Whitehead link
in \cite{L}.
 We briefly recall his tactics here.
 We calculate boundary slopes of surfaces
with respect to the meridian $\mu_i$ and a non-standard longitude $\lambda_i$
of $L_i$.
 We take $\lambda_1$ as below.
 In $S^2 = {\Bbb R}^2 / \Gamma$,
we take the arc $s$ of slope $0$ connecting $\Gamma(0,0)$ and $\Gamma(0,1)$.
 $\lambda_1$ is the union of the arc $(s \times [0,1]) \cap \partial N(L_1)$
and an arc in $(S^2 \times [1, \infty)) \cap \partial N(L_1)$.
 $\lambda_1$ is oriented
toward increasing $r \in [0,1]$ along the axis $\Gamma(0,0) \times [0,1]$,
and hence toward decreasing along the axis $\Gamma(0,1) \times [0,1]$.
 The meridian $\mu_1$ is oriented
as a right-handed circle around the axis $\Gamma(0,0) \times [0,1]$
oriented upward.
 We obtain $\lambda_2$, $\mu_2$ from $\lambda_1$, $\mu_1$
by rotating by $180^{\circ}$ about the axis $\Gamma(1/2,1/2) \times [0,1]$.
 (Under this coordination,
the preferred longitude of $L([r,s])$ is
of slope $(1, \frac{r-s}{2})$ for both components when $s>0$,
as is shown in section 7.
 It is $(1, \frac{r-s}{2} -2)$ when $s<0$, as is shown in section 8.)

 Let $i_j$ be the algebraic intersection number $\partial F \cdot \lambda_j$
in $\partial N(L_i)$.
 For $1 < t \le \infty$,
$\partial \Sigma_{e_i} = \partial (g \times *) (\Sigma_{e_0})$,
$g = \left( \begin{array}{cc}
a & b \\
c & d \\
\end{array} \right) \in G$
contributes to the number $i_j$ as in Table 2.1
(Table 1.2 in \cite{L})
if the orientations of $e_i$ and $g(e_0)$ agree.
 If they disagree,
we change all the sign of the number in Table 2.1.
 (If $\lambda_i$ is not transverse to $\partial \Sigma_{e_i}$ for some $i$,
we isotope whole of $\lambda_i$ very slightly to the right
in $\partial N(L_i)$.)
 Fortunately, Table 2.1 does not depend on the branching number $n$.
 For $t=\infty$, we substitute $0$ for $\beta$ in Table 2.1.

\vspace{3mm}
\begin{center}
$\begin{array}{|c|c|c|c|}
\hline
{\rm label} & {\rm condition\ on}\ g & i_1 & i_2 \\
\hline
A & -\infty < -\frac{d}{c} < 0 &  \beta & \beta \\
\cline{2-4}
  & 0 < -\frac{d}{c} < \infty & -\beta & -\beta \\
\cline{2-4}
  & -\frac{d}{c} = 0, \pm \infty & 0 & 0 \\
\hline
B & -\infty < -\frac{d}{c} < 0 & -(\alpha-\beta) & 0 \\
\cline{2-4}
  & 0 < -\frac{d}{c} < \infty & \alpha-\beta & 0 \\
\cline{2-4}
  & -\frac{d}{c} = 0, \pm \infty & 0 & 0 \\
\hline
C & 0 < -\frac{d}{c} < 1 & -2\beta & 0 \\
\cline{2-4}
  & -\frac{d}{c} = 0,1 & -\beta & \beta \\
\cline{2-4}
  & otherwise & 0 & 2\beta \\
\hline
D & \frac{1}{2} < -\frac{d}{c} < \infty & \alpha-\beta & \alpha-\beta \\
\cline{2-4}
  & -\frac{d}{c} = \frac{1}{2}, \pm \infty & 0 & \alpha-\beta \\
\cline{2-4}
  & {\rm otherwise} &  -(\alpha-\beta) & \alpha-\beta \\
\hline
\end{array}$\\
Table 2.1
\end{center}

 For $0 \le t < 1$,
we interchange the column for $i_1$ and that for $i_2$
and interchange $\alpha$ and $\beta$ in Table 2.1.
 For $t=0$, we further substitute $0$ for $\alpha$.

 For $t=1$,
contribution of $\Sigma_{e_i} = (g \times *) \Sigma_A$ to $i_1, i_2$
is the same as Table 2.1.
 $\Sigma_C$ corresponds to the edge from $1/1$ to $0/1$,
and 
$\Sigma_C
=
(\left( \begin{array}{cc}
1 & 0 \\
1 & 1 \\
\end{array} \right)
\times id) \Sigma_A$.
 Hence, for an edge $e_i$ labeled by $C$,
$\partial \Sigma_{e_i} = \partial g (\Sigma_C)$
with
$g =
\left( \begin{array}{cc}
a & b \\
c & d \\
\end{array} \right)
\in G$
contributes $i_1$ and $i_2$ as in Table 2.2
(Table 1.3 in \cite{L}).
 We should note that it depends on $n$.

\vspace{3mm}
\begin{center}
$\begin{array}{|c|c|c|c|}
\hline
  & {\rm condition\ on}\ g & i_1 & i_2 \\
\hline
C & -\frac{d}{c} = 0,1 & \beta-2n & -(\beta-2n) \\
\cline{2-4}
  & 0 < -\frac{d}{c} < 1 & -2n & -2(\beta-n) \\
\cline{2-4}
  & {\rm otherwise} & 2(\beta-n) & 2n \\
\hline
\end{array}$\\
Table 2.2
\end{center}

\vspace{3mm}

 In the remainder of this section,
we recall the way how we obtain Tables 2.1 and 2.2.
 We consider, for example,
the case where $1 < t < \infty$ and $e_i = g (D_0)$, $g \in G$
with orientations match.
 Set $g = \left( \begin{array}{cc} a & b \\ c & d \\ \end{array} \right)$.
 Then $g$ brings the slope $-\dfrac{d}{c}$ to $1/0$
since
$\left( \begin{array}{cc} a & b \\ c & d \\ \end{array} \right)
\left( \begin{array}{c} -d \\ c \\ \end{array} \right)
=\left( \begin{array}{c} -ad+bc \\ 0 \end{array} \right)
=\left( \begin{array}{c} -1 \\ 0 \end{array} \right)$.
 Hence the intersection number between $\partial \Sigma_{e_i}$ and $\lambda_j$
is equal to that between $\partial \Sigma_D$
and the vertical line segment $(\rho \times [0,1]) \cap \partial N(L_j)$
where $\rho$ is the arc of slope $-\dfrac{d}{c}$
in $S^2 = {\Bbb R}^2 / \Gamma$.

\begin{figure}[htbp]
\begin{center}
\includegraphics[width=9cm]{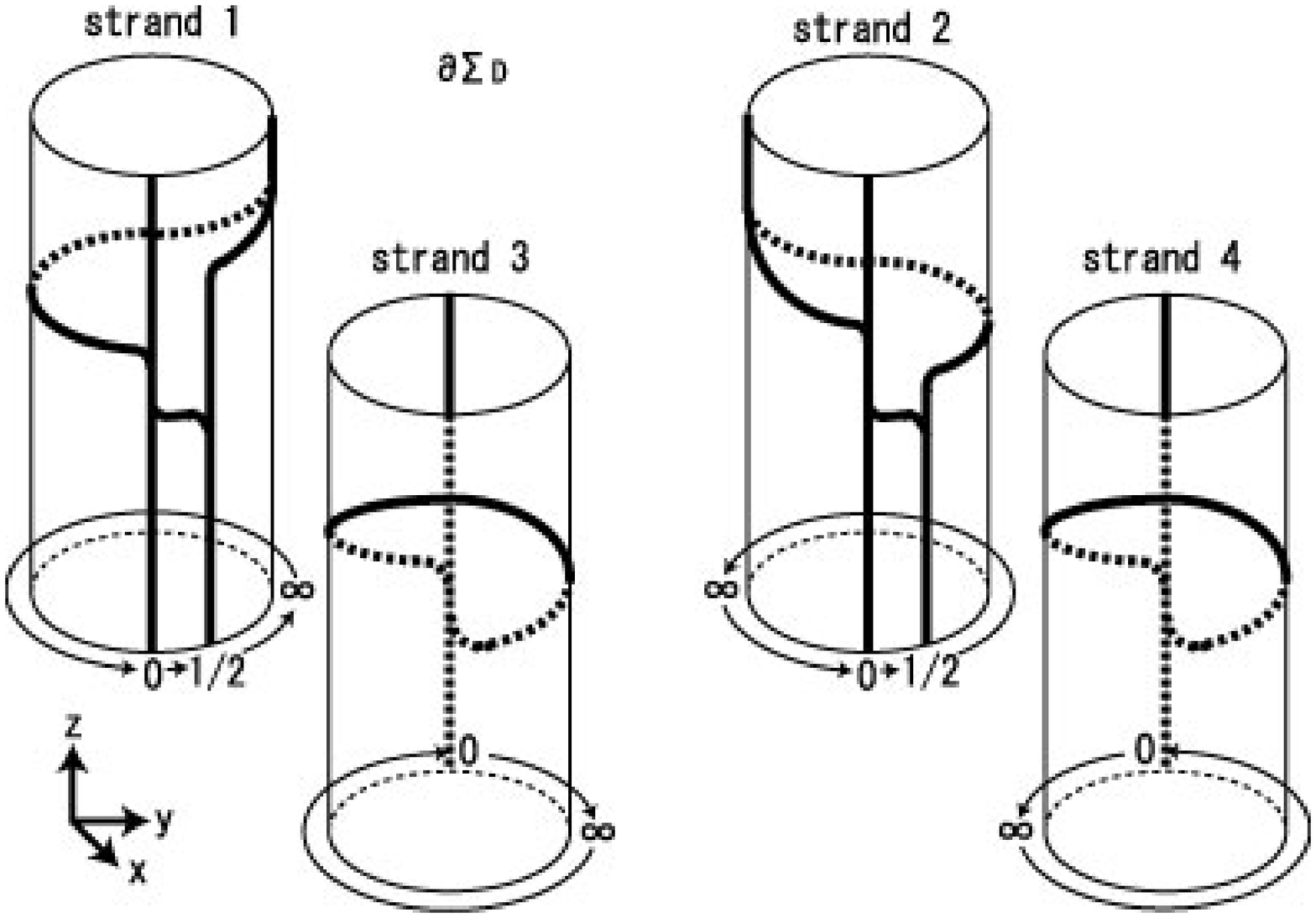}
\end{center}
\caption{}
\label{fig:SigmaD1}
\end{figure}

\begin{figure}[htbp]
\begin{center}
\includegraphics[width=10cm]{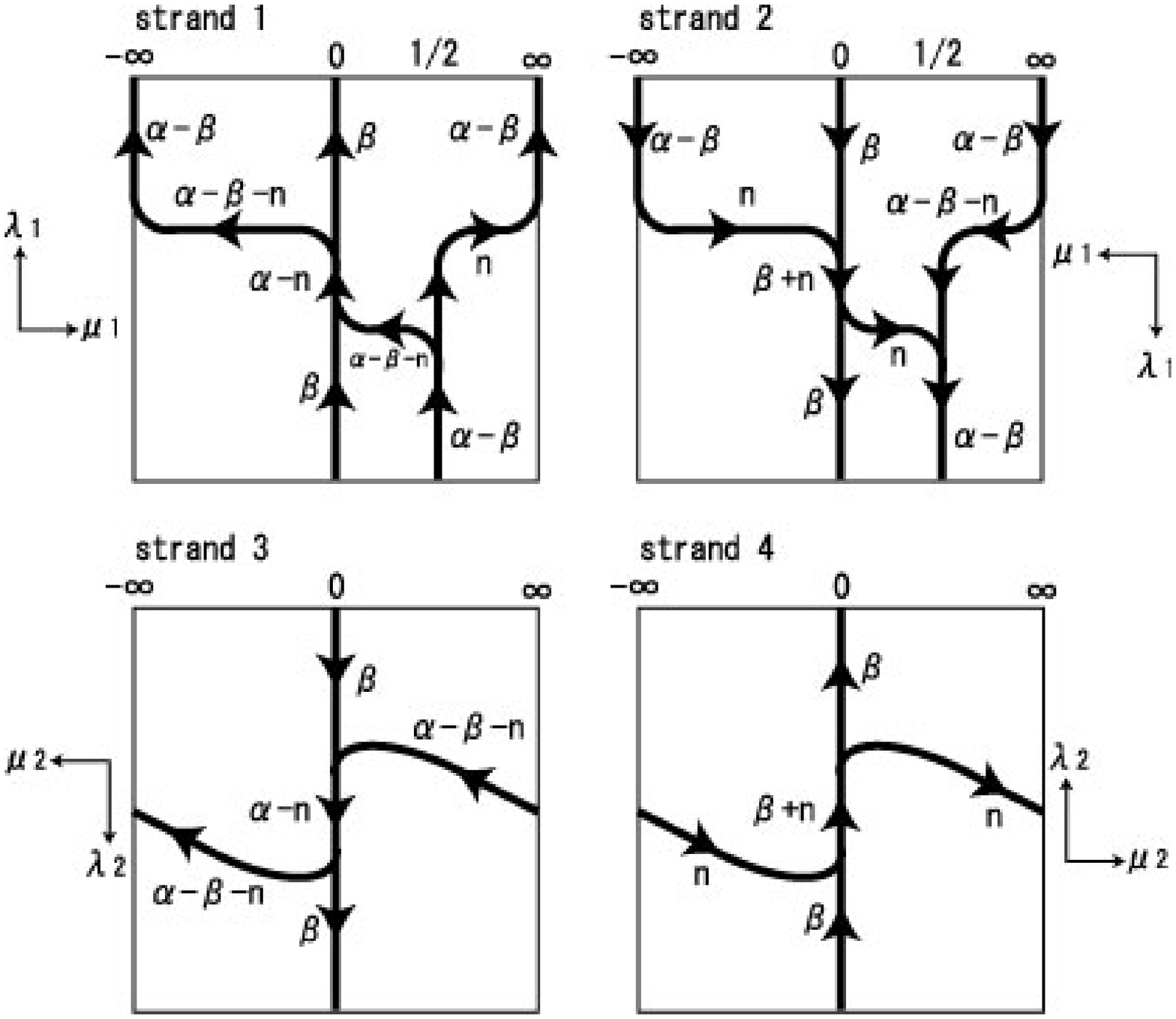}
\end{center}
\caption{}
\label{fig:SigmaD2}
\end{figure}

 In Figure \ref{fig:SigmaD1},
we can find $\partial \Sigma_D$.
 We orient $\partial \Sigma_{\gamma}$
so that it runs
along the strands $\Gamma(0,0) \times [0,1]$ and $\Gamma(1,1) \times [0,1]$
to point upward.
 Figure \ref{fig:SigmaD2} depicts the tubes about the four strands
lifted into ${\Bbb R}^2 \times [0,1]$ via $\Gamma$.

 We take the algebraic sum of these signed crossing points,
to obtain Table 2.3 below.
 Note that $\lambda_j$ points downward
along strand 2 $\Gamma(0,1) \times [0,1]$
and strand 3 $\Gamma(1,0) \times [0,1]$.
 We should recall
that, if $\lambda_j$ is not transverse to $\partial \Sigma_{e_i}$,
then $\lambda_j$ is isotoped slightly to the direction of orientation of $\mu_j$
if $\lambda_j$ is not transverse to $\partial \Sigma_{e_i}$.
 Then we take sum of the intersection numbers on strands 1 and 2,
and those on strands 3 and 4, to obtain the bottom three rows in Table 2.1.

\vspace{3mm}
\begin{center}
$\begin{array}{|c|c|c|c|c|c|}
\hline
{\rm label} & {\rm condition on} g & i_1, {\rm strand 1} & i_1, {\rm strand 2}
& i_2, {\rm strand 3} & i_2, {\rm strand 4} \\
\hline
D & \frac{1}{2} < -\frac{d}{c} < \infty & n & \alpha-\beta-n
& \alpha-\beta-n & n \\
\cline{2-6}
 & -\frac{d}{c}=\frac{1}{2} & n & -n
& \alpha-\beta-n & n \\
\cline{2-6}
 & -\frac{d}{c} = \infty & -(\alpha-\beta-n) & \alpha-\beta-n
& \alpha-\beta-n & n \\
\cline{2-6}
 & {\rm otherwise} & -(\alpha-\beta-n) & -n
& \alpha-\beta-n & n \\
\hline
\end{array}$ \\
Table 2.3
\end{center}

\section{Minimal edge paths}

 In this section we list all the minimal edge-paths for $L(p/q)$
with $p/q=[r,s]=1/(r+1/s)$, $r = 2w+1$, $s = 2u+1$, $r \ge 3$ and $s \ne \pm 1$.
 Since $L(p/q)=L([r,s])$ and $L(-p/q)=L([-r,-s])$
are mirror images of each other,
we can assume $r > 0$.

\begin{figure}[htbp]
\begin{center}
\includegraphics[width=8cm]{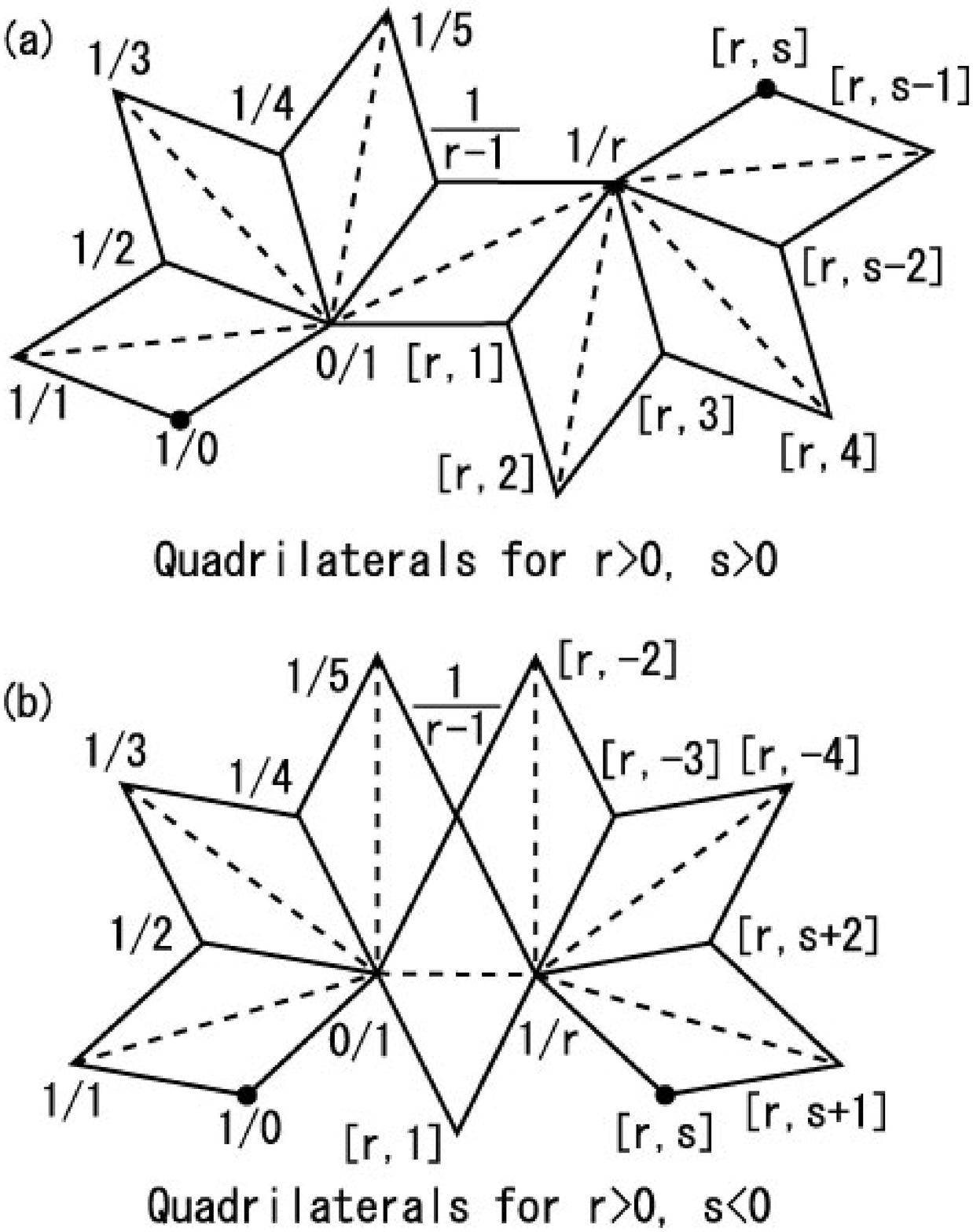}
\end{center}
\caption{}
\label{fig:quadra}
\end{figure}

\noindent
{\bf The case of $s>0$.}

 We consider the case of $u \ge 1$.
 Then $Q_{p/q}$ is a union of
the sequence of $w+1$ quadrilaterals around the vertex $0/1$,
and the sequence of $u+1$ quadrilaterals around the vertex $1/r$.
 See Figure \ref{fig:quadra} (a).
 The first sequence is composed of the quadrilaterals
$\langle  0/1, 1/0, 1/1, 1/2 \rangle $,
$\langle 0/1, 1/2, 1/3, 1/4 \rangle $, $\cdots$,
$\langle 0/1, 1/(2i-2), 1/(2i-1), 1/2i \rangle $, $\cdots$,
$\langle 0/1, 1/(r-1), 1/r, 1/(r+1) \rangle $
 The last sequence is the union of the quadrilaterals
$\langle 1/r, [r,-1], [r,0], [r,1] \rangle $,
$\langle 1/r, [r,1], [r,2], [r,3] \rangle $, $\cdots$,
$\langle 1/r, [r,2i-1], [r,2i], [r,2i+1] \rangle $, $\cdots$,
$\langle 1/r, [r,s-2], [r,s-1], [r,s] \rangle $.
 Since $[r,0]=1/(r+1/0)=1/\infty=0=0/1$,
the last quadrilateral of the first sequence
is the first one of the last sequence.
 Hence $Q_{p/q}$ consists of $w+u+1$ quadrilaterals.

\begin{figure}[htbp]
\begin{center}
\includegraphics[width=12cm]{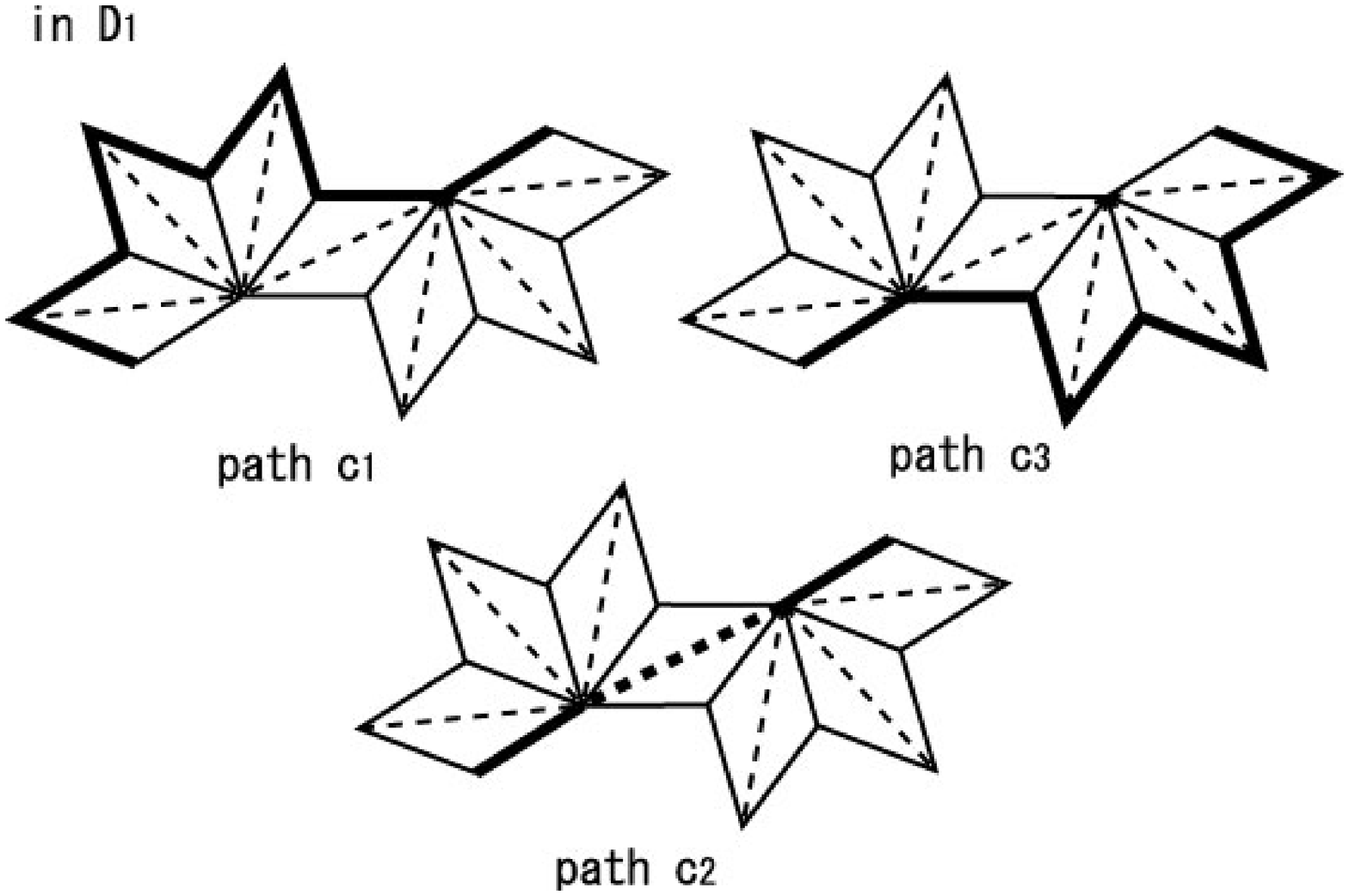}
\end{center}
\caption{}
\label{fig:PathPD1}
\end{figure}

 For the diagram $D_1$,
we have $3$ minimal edge paths.
 See Figure \ref{fig:PathPD1}.

 The edge-path $c_1$ is a sequence of $2w+2$ $A$-edges.
 The first $2w$ edges connect the vertices
$1/0$-$1/1$-$1/2$-$1/3$-...-$1/(r-1)$,
and the last $2$ edges connect $1/(r-1)$-$1/r$-$[r,s]$.

 The edge-path $c_2$ is a sequence of $3$ edges.
 The first $A$-edge connects $1/0$ and $0/1$,
the second $C$-edge connects $0/1$ and $1/r$
and the last $A$-edge connects $1/r$ and $[r,s]$.

 The edge-path $c_3$ is a sequence of $2u+2$ $A$-edges.
 $\pi$-rotation of $c_3$ on the paper is $c_1$ in $Q_{[s,r]}$.

\begin{figure}[htbp]
\begin{center}
\includegraphics[width=12cm]{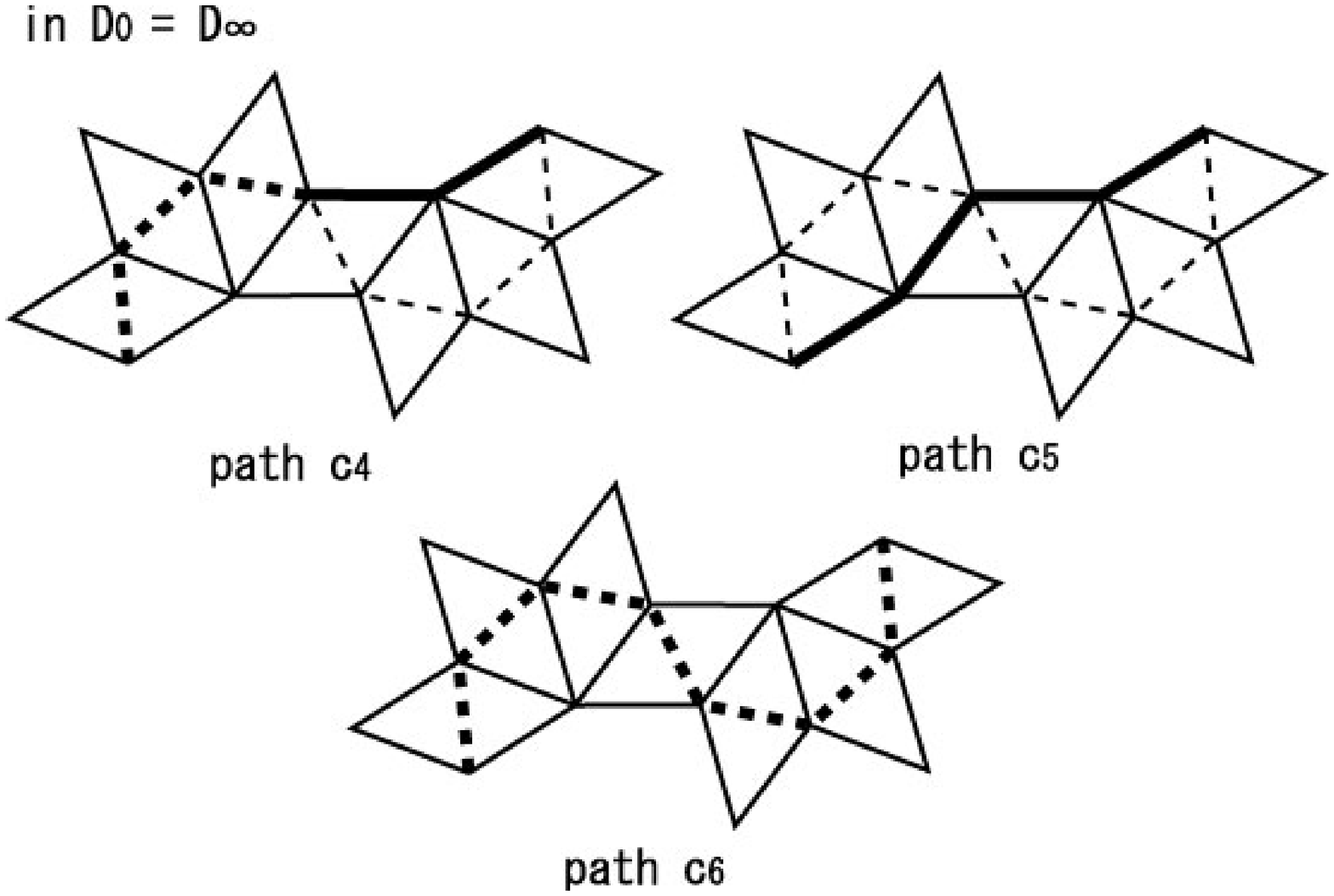}
\end{center}
\caption{}
\label{fig:PathPD0}
\end{figure}

 For the diagram $D_0 = D_{\infty}$,
we have $5$ minimal edge-paths.
 See Figure \ref{fig:PathPD0}.

 The edge-path $c_4$ is a union of $2$ sequence of edges.
 The first sequence of $w$ $D$-edges
connect $1/0$-$1/2$-$1/4$-$1/6$-...-$1/(2w)=1/(r-1)$.
 The second is a sequence of $2$ $B$-edges
connecting $1/(r-1)$-$1/r$-$[r,s]$.

 The edge-path $c_5$ is the union of $4$ $B$-edges
connecting $1/0$-$0/1$-$1/(r-1)$-$1/r$-$[r,s]$.
 This edge-path is not minimal when $r=3$ ($w=1$).

 The edge-path $c_6$ is a union of $w+u+1$ $D$-edges.
 The first $t$ edges connect $1/0$-$1/2$-$1/4$-$1/6$-...-$1/(2w)=1/(r-1)$.
 The middle edge connect $1/(r-1)$-$1/(r+1)$.
 The last $u$ edges connect $1/(r+1)$-$[r,3]$-...-$[r,s]$.

 The edge-paths $c_7$ and $c_8$ correspond
to the edge-paths $c_5$ and $c_4$ respectively
by the $\pi$-rotation on the paper.

\begin{figure}[htbp]
\begin{center}
\includegraphics[width=12cm]{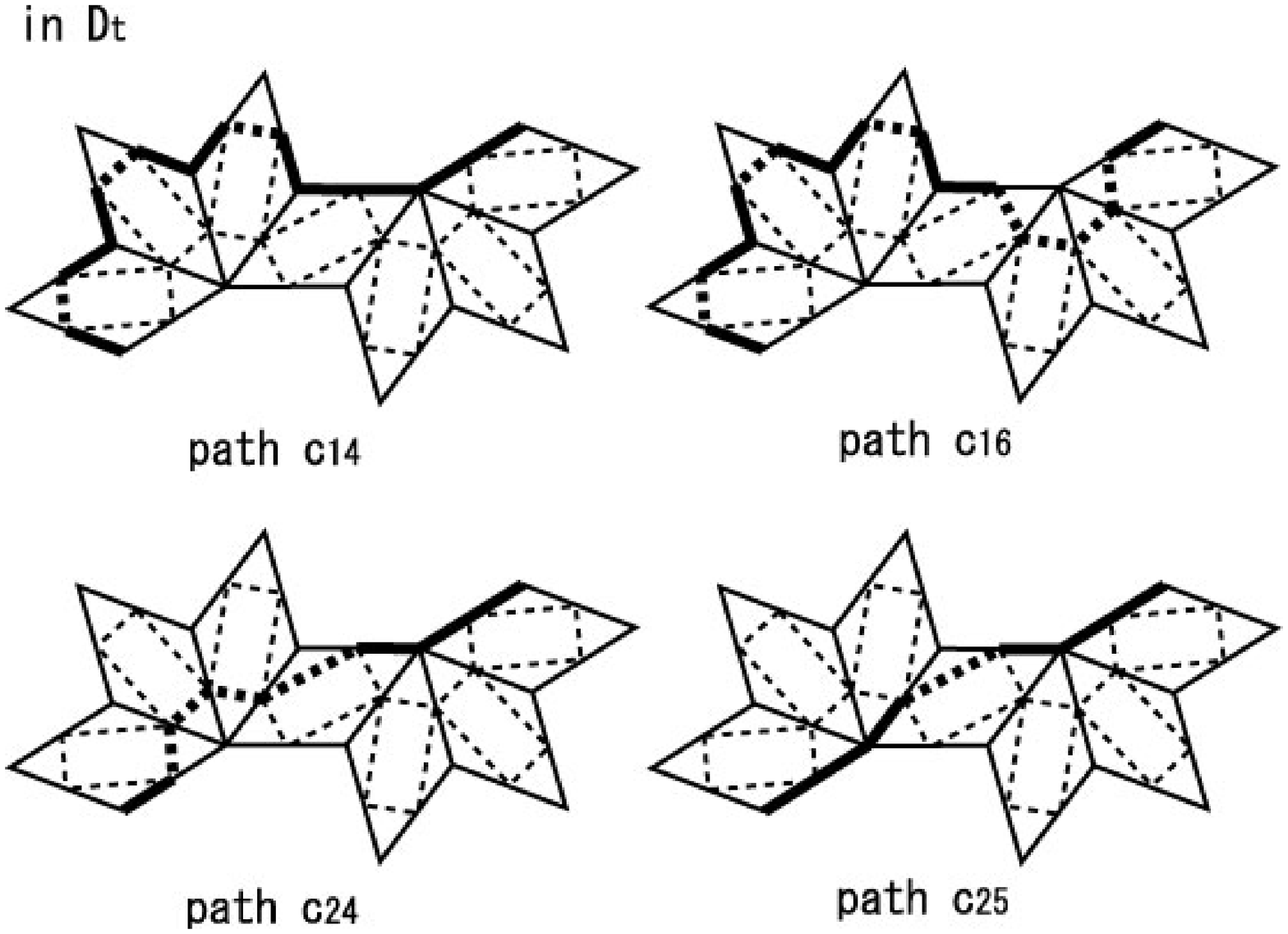}
\end{center}
\caption{}
\label{fig:PathPDt}
\end{figure}

 For the diagram $D_t$ with $t \ne 0,1,\infty$,
we have $8$ minimal edge-paths.
 See Figure \ref{fig:PathPDt}.
 The edge-path $c_{ij}$ is a $\lq\lq$chimera" of $c_i$ and $c_j$.

 The edge-path $c_{14}$ is composed of $2$ sequences of edges.
 The first is the sequence of $w$ triples of $A$-, $D$-, $A$-edges
together connecting the vertices $1/0$-$1/2$-$1/4$-...-$1/(2w)=1/(r-1)$.
 These edges are away from the vertex $0/1$.
 The last is the sequence of $4$ $A$-, $B$-, $B$-, $A$-edges,
connecting $1/(r-1)$-$1/r$-$1/[r,s]$.

 The edge-path $c_{16}$ is a union of $2$ sequences of edges.
 The first is the same as that of $c_{14}$.
 The last is the sequence of an $A$-edge, $u+1$ $D$-edges and an $A$-edge
together connecting $1/(r-1)$ to $[r,s]$ near the vertex $1/r$.

 The edge-path $c_{24}$ is a union of $2$ sequences of edges.
 The first is the sequence of $A$-, $w$ $D$-, $C$-, $B$-edges
connecting $1/0$-$1/r$,
where the $D$-edges are near the vertex $0/1$.
 The last consists of a $B$-edge and an $A$-edge connecting $1/r$ to $[r,s]$.

 The edge-path $c_{25}$ is the union of 7 edges
which are $A$-, $B$-, $B$-, $C$-, $B$-, $B$-, $A$-edges
connecting $1/0$-$0/1$-$1/r$-$[r,s]$.
 The $C$-edge is contained in the quadrilateral $\langle 0/1, 1/(r-1), 1/r, 1/(r+1) \rangle $
and near the vertex $1/(r-1)$ rather than $1/(r+1)$.
 This edge-path is not minimal when $r=3$ ($t=1$).

 The edge-paths $c_{38}$, $c_{36}$, $c_{28}$, $c_{27}$ correspond
to the edge-paths $c_{14}$, $c_{16}$, $c_{24}$, $c_{25}$ respectively
by the $\pi$-rotation.


\vspace{3mm}

\noindent
{\bf The case of $s<0$.}

 We set $s=2u+1=-(2u'+1)$.
 Then $u \le -2$ and $u' \ge 1$.
 The sequence of quadrilaterals $Q_{p/q}$
is the union of the 2 sequences of quadrilaterals below.
 See Figure \ref{fig:quadra} (b).
 The first sequence is the precisely same one as that for $s>0$,
and is composed of $w+1$ quadrilateral.
 The last sequence consists of $u'+1$ quadrilaterals around the vertex $1/r$.
 Precisely, it is the union of the quadrilaterals
$\langle 1/r, [r,1], [r,0], [r,-1] \rangle $,
$\langle 1/r, [r,-1], [r,-2], [r,-3] \rangle $, $\cdots$
$\langle 1/r, [r,-2i+3], [r,-2i+2], [r,-2i+1] \rangle $, $\cdots$
$\langle 1/r, [r,s+2], [r,s+1], [r,s] \rangle $.
 Thus $Q_{p/q}$ consists of $w+u'+1=w-u$ quadrilaterals.

\begin{figure}[htbp]
\begin{center}
\includegraphics[width=12cm]{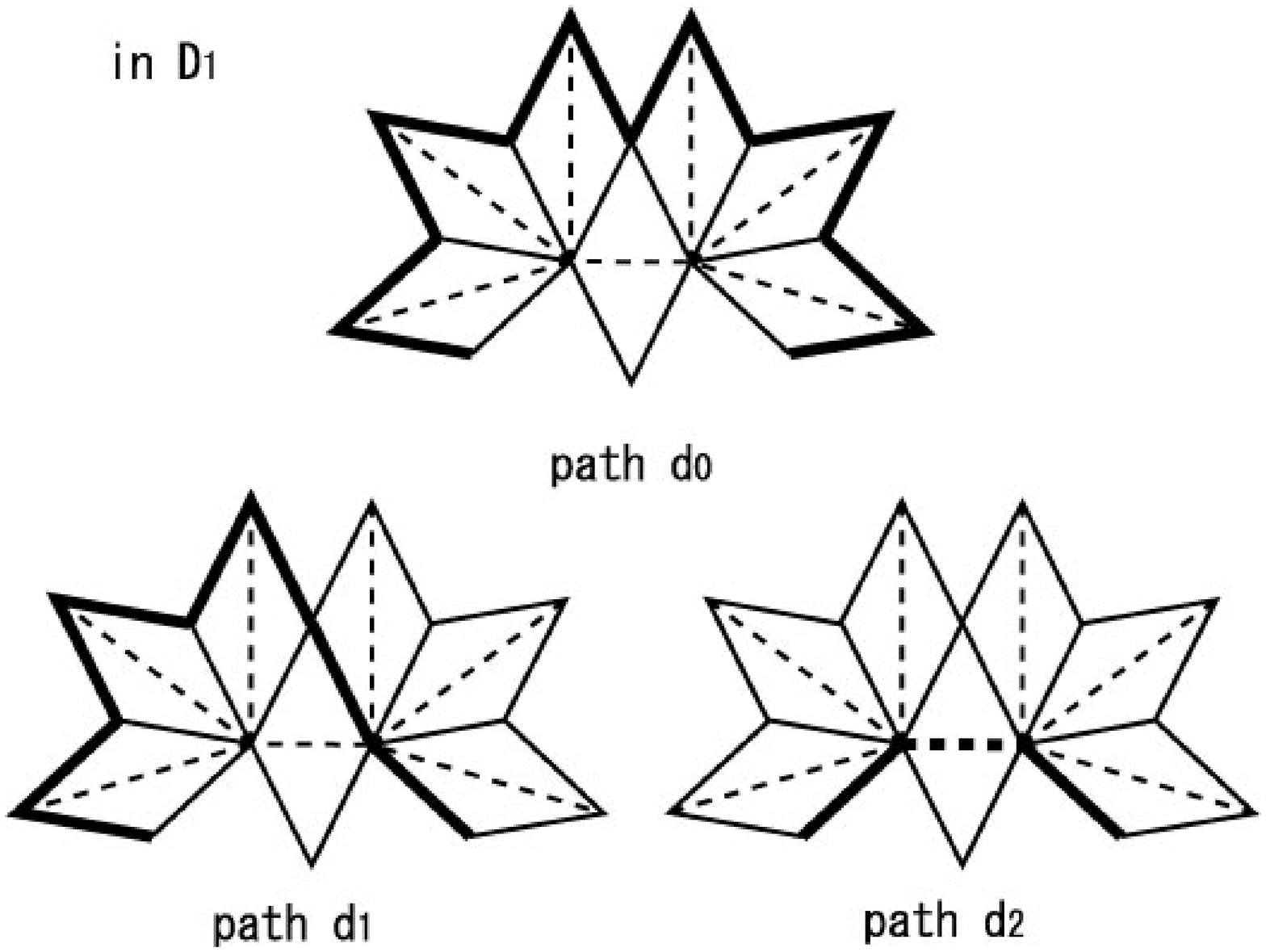}
\end{center}
\caption{}
\label{fig:PathND1}
\end{figure}

 For the diagram $D_1$,
we have $4$ minimal edge-paths.
 See Figure \ref{fig:PathND1}.

 The edge-path $d_0$ is a union of $2$ sequences of edges.
 The first is the sequence of $2w$ $A$-edges
connecting $1/0$-$1/1$-$1/2$-$1/3$-$1/4$-$1/5$-...-$1/(r-1)$.
 The last is the sequence of $2u'$ $A$-edges
connecting $1/(r-1)=[r,-1]$-$[r,-2]$-$[r,-3]$-...-$[r,s]$.
 Thus the edge-path consists of $2w+2u'$ edges.

 The edge-paths $d_1$, $d_2$ and $d_3$ are similar
to $c_1, c_2$ and $c_3$ for $s>0$.
 The edge-path $d_3$ corresponds to the edge-path $d_1$ for $Q_{[-s,-r]}$
by reflection on a vertical axis on the paper.

\begin{figure}[htbp]
\begin{center}
\includegraphics[width=12cm]{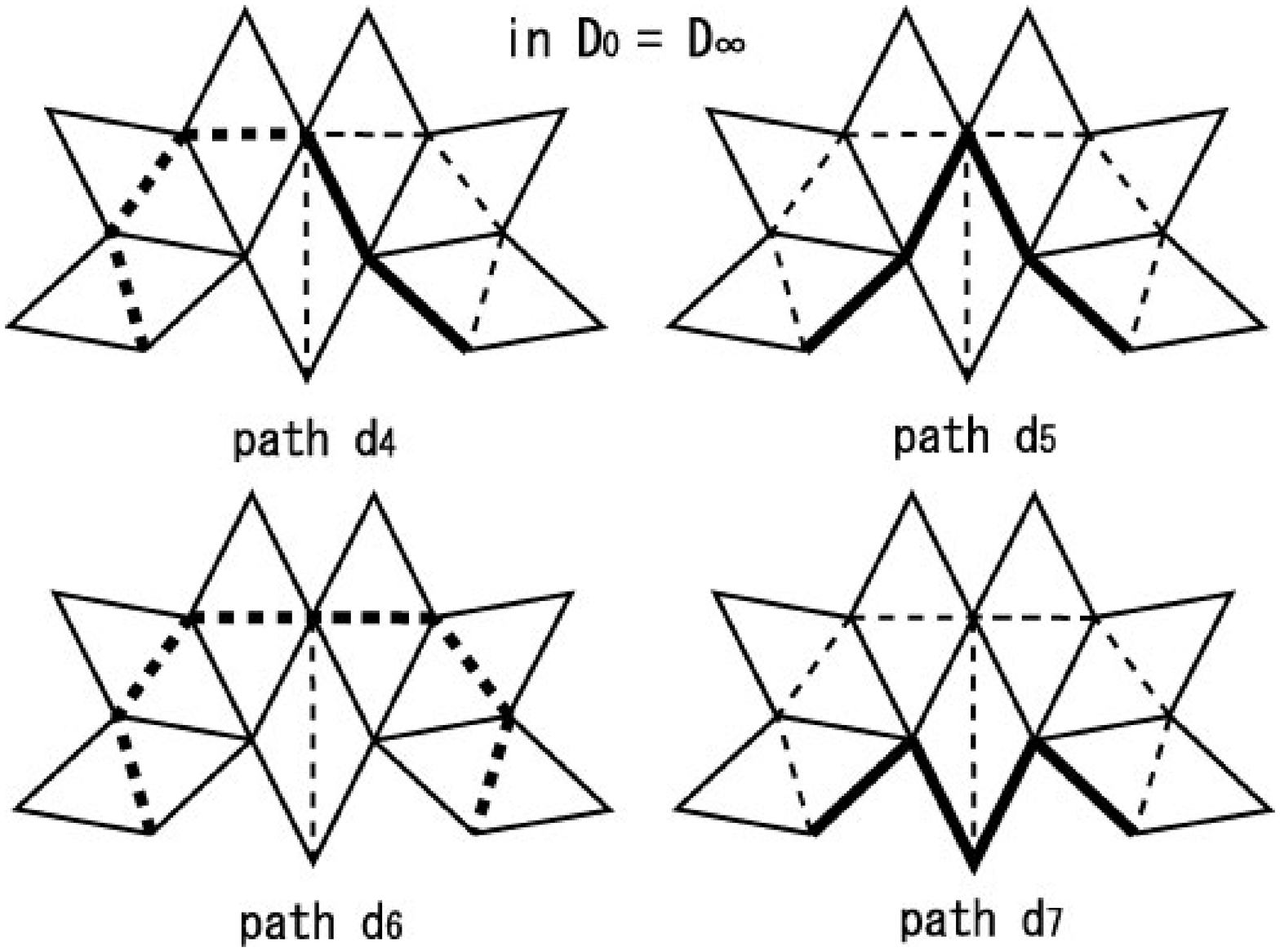}
\end{center}
\caption{}
\label{fig:PathND0}
\end{figure}

 For the diagram $D_0 = D_{\infty}$,
we have $5$ minimal edge-paths.
 See Figure \ref{fig:PathND0}.

 The edge-paths $d_4$, $d_5$, $d_6$, $d_7$ and $d_8$ are similar
to $c_4, c_5, c_6, c_7$ and $c_8$ for $s>0$.
 But, when $r=3$ ($w=1$), $d_5$ and $d_8$ are not minimal.
 $d_4$ and $d_5$ are not minimal when $s=-3$ ($u'=1$, $u=-2$).
 (Note that $c_7$ is minimal even when $r=3$ or $s=-3$.)
 By reflection on a vertical axis,
$d_8$ corresponds to $d_4$ for $Q_{[-s,-r]}$,
while $d_7$ does not correspond to $d_5$.

\begin{figure}[htbp]
\begin{center}
\includegraphics[width=12cm]{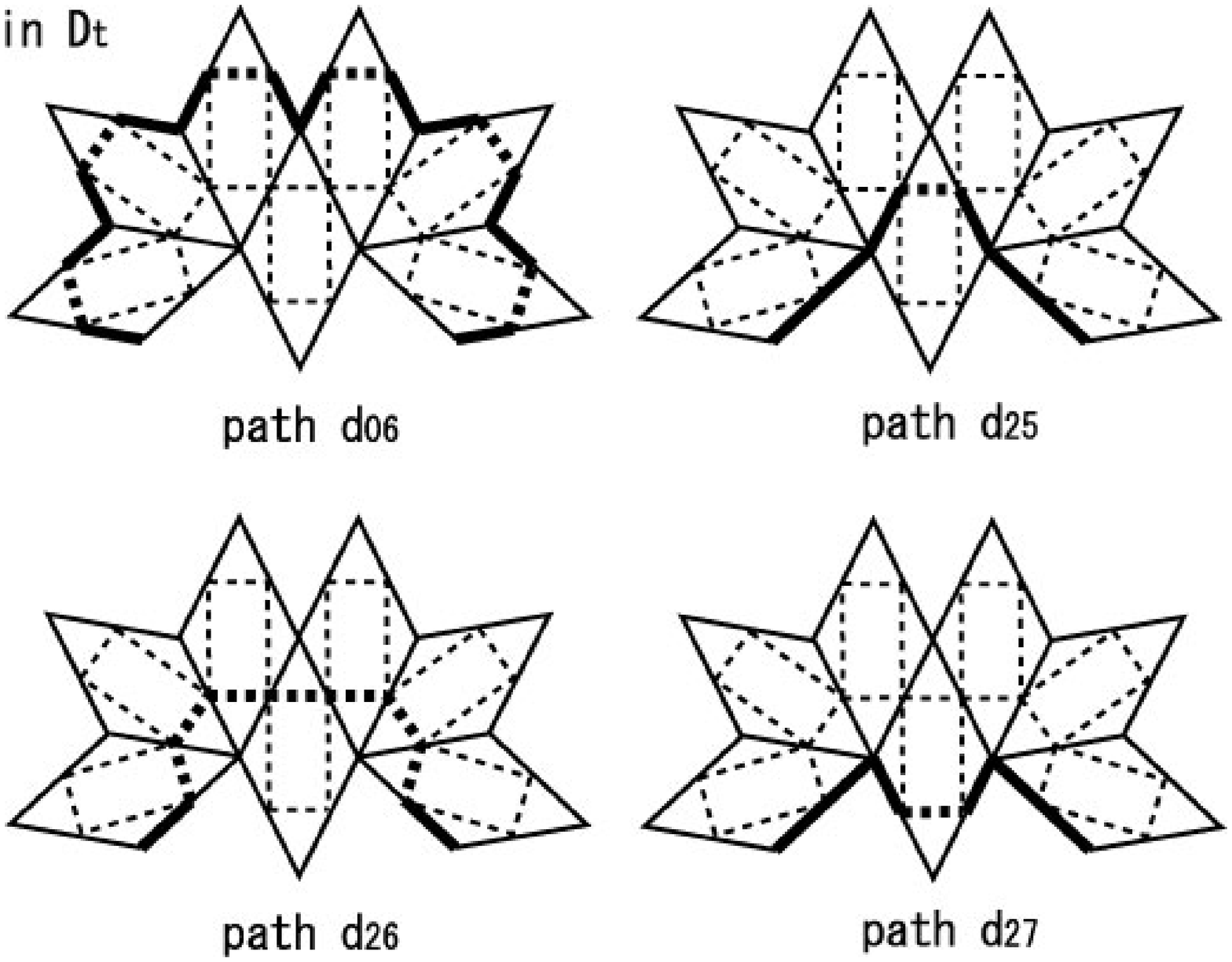}
\end{center}
\caption{}
\label{fig:PathNDt0}
\end{figure}

\begin{figure}[htbp]
\begin{center}
\includegraphics[width=12cm]{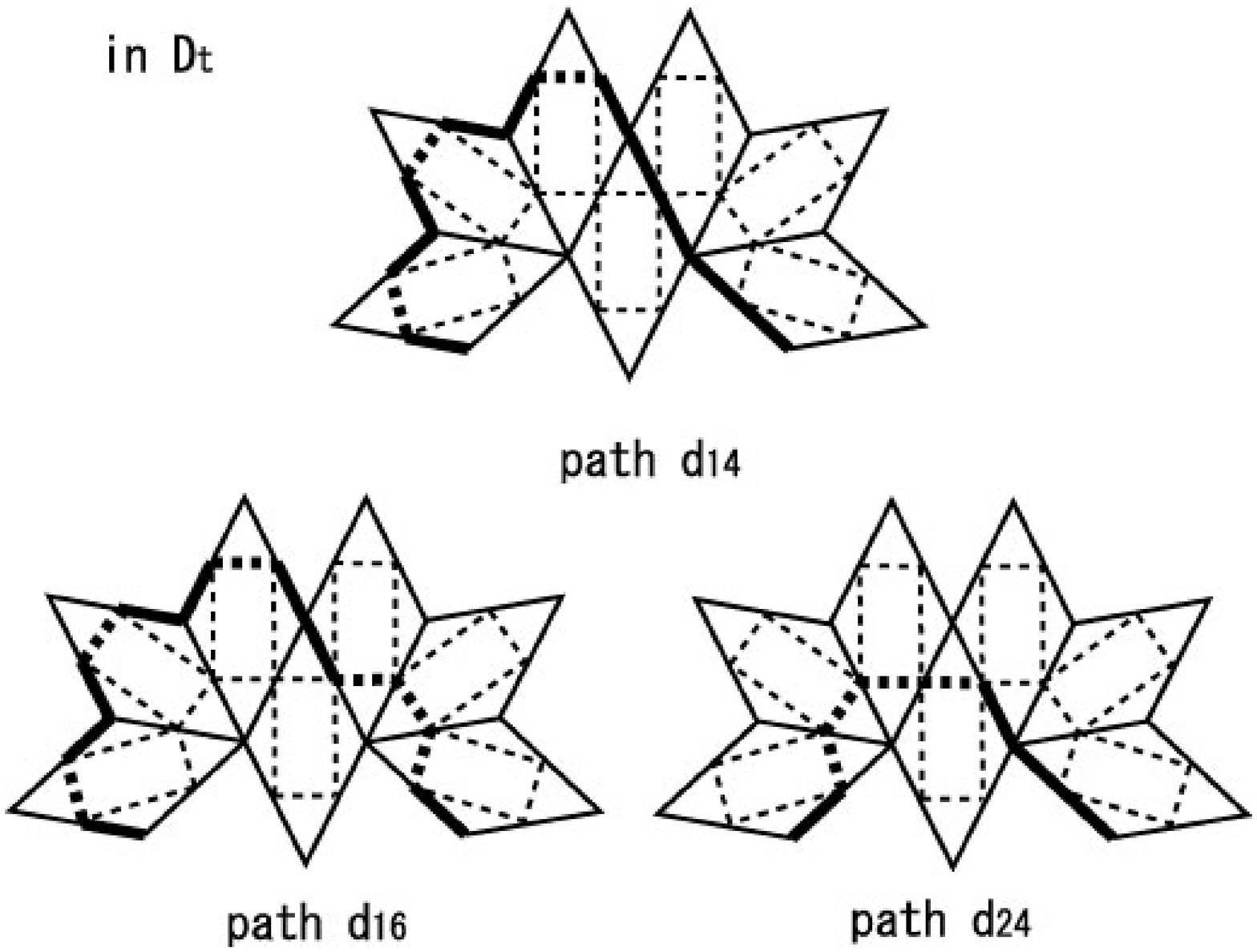}
\end{center}
\caption{}
\label{fig:PathNDt1}
\end{figure}

 For the diagram $D_t$ with $t \ne 0,1,\infty$,
we have $10$ minimal edge paths.
 See Figures \ref{fig:PathNDt0} and \ref{fig:PathNDt1}.

 The edge-path $d_{06}$ is a union of $2$ sequences of edges.
 The first is the sequence of $w$ $A$-, $D$-, $A$-edges
connecting $1/0$-$1/2$-$1/4$-...-$1/(2w)=1/(r-1)$
away from the vertex $0/1$.
 The last is the sequence of $u'$ $A$-, $D$-, $A$-edges
connecting $[r,-1]$-$[r,-3]$-$[r,-5]$-...-$[r,s]$
away from the vertex $1/r$.
 Hence this edge-path consists of $3(w+u')$ edges.

 The edge-path $d_{26}$ is a sequence of
an $A$-edge, $w$ $D$-edges, a $C$-edge, $u'$ $D$-edges and an $A$-edge.
 The first sequence of $w$ $D$-edges are near the vertex $0/1$,
and the last ones are near $1/r$.

 For the other $8$ edge-paths
$d_{14}$, $d_{16}$, $d_{24}$, $d_{25}$, $d_{38}$, $d_{36}$, $d_{28}$, $d_{27}$,
each $d_{ij}$ is similar to $c_{ij}$ for $s>0$.
 $d_{ij}$ is minimal if and only if $d_j$ is minimal.


\section{transformation $g \in G$ with $e_i = g(e_0)$ and intersection number}

\begin{figure}[htbp]
\begin{center}
\includegraphics[height=20cm]{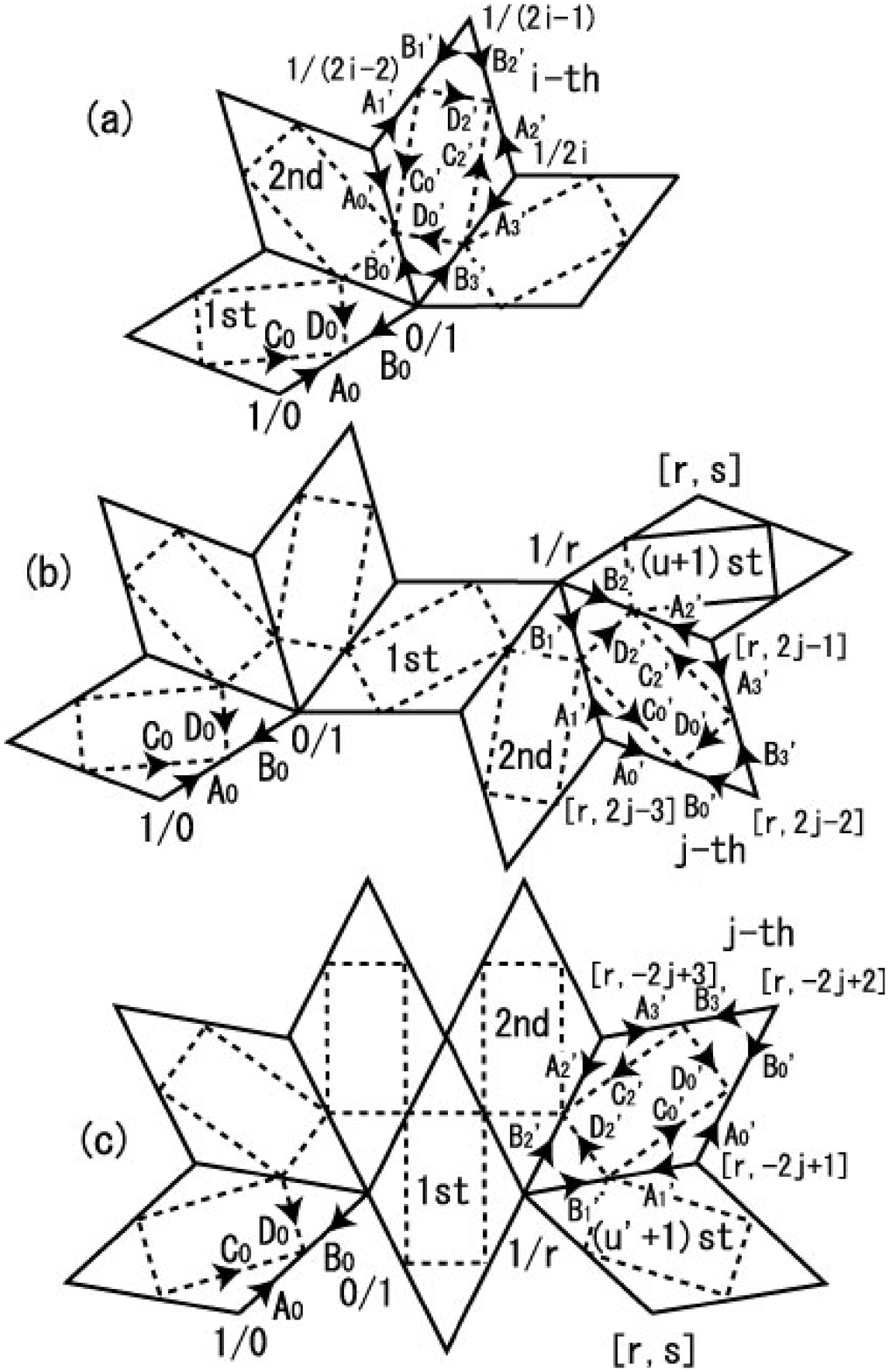}
\end{center}
\caption{}
\label{fig:trans}
\end{figure}

 For $E \in \{ A, B, C, D \}$,
the edges $E'_i$ in Figure \ref{fig:trans} are the images of $E_0$
under the following transformations.

 For the $i$-th quadrilateral in the first sequence of quadrilaterals $Q_{p/q}$,

$E'_0 = f_0 (E_0)$,
$f_0 =
\left( \begin{array}{cc}
1 & 0 \\
2i-2 & 1 \\
\end{array} \right)
\in G$,
$-\dfrac{d}{c} = -\dfrac{1}{2i-2}$,

$E'_1 = f_1 (E_0)$,
$f_1 =
\left( \begin{array}{cc}
1 & 1 \\
2i-2 & 2i-1 \\
\end{array} \right)
\in G$,
$-\dfrac{d}{c} = -\dfrac{2i-1}{2i-2}$,

$E'_2 = f_2 (E_0)$,
$f_2 =
\left( \begin{array}{cc}
1 & -1 \\
2i & -(2i-1) \\
\end{array} \right)
\in G$,
$-\dfrac{d}{c} = \dfrac{2i-1}{2i}$, and

$E'_3 = f_3 (E_0)$,
$f_3 =
\left( \begin{array}{cc}
1 & 0 \\
2i & 1 \\
\end{array} \right)
\in G$,
$-\dfrac{d}{c} = -\dfrac{1}{2i}$

 The first column $\left( \begin{array}{c} a_k \\ c_k \end{array} \right)$
of $f_k$
corresponds to the initial vertex $a_k / c_k$ of the edge $A_k$,
and the second column $\left( \begin{array}{c} b_k \\ d_k \end{array} \right)$
of $f_k$
corresponds to the initial vertex $b_k / d_k$ of the edge $B_k$.
 The sign of each column is determined
so that the determinant of $f_k$ is equal to $+1$ rather than $-1$.
 Note that $c_k$ is an even integer.

 Hence, when $1 \le t \le \infty$ ($1 < t \le \infty$ for $C'_k$),
from Table 2.1,
we obtain Table 4.1 of contribution
to the intersection numbers $i_1$ and $i_2$
for the edges of the $i$-th quadrilateral
in the first sequence of $Q_{p/q}$.

\begin{center}
$\begin{array}{|c|c|c|c|c|c|c|}
\hline
 & \multicolumn{3}{c|}{i=1} & \multicolumn{3}{c|}{i \ge 2} \\
\hline
 & -\frac{d}{c} & i_1 & i_2 & -\frac{d}{c} & i_1 & i_2 \\
\hline
A'_0 & \infty & 0 & 0 & - & \beta & \beta \\
\hline
B'_0 & \infty & 0 & 0 & - & -(\alpha-\beta) & 0 \\
\hline
C'_0 & \infty & 0 & 2\beta & - & 0 & 2\beta \\
\hline
D'_0 & \infty & 0 & \alpha-\beta & - & -(\alpha-\beta) & \alpha-\beta \\
\hline
A'_1 & \infty & 0 & 0 & - & \beta & \beta \\
\hline
B'_1 & \infty & 0 & 0 & - & -(\alpha-\beta) & 0 \\
\hline
A'_2 & \frac{1}{2} & -\beta & -\beta 
& \frac{1}{2}< -\frac{d}{c} < 1 & -\beta & -\beta \\
\hline
B'_2 & \frac{1}{2} & \alpha-\beta & 0
& \frac{1}{2}< -\frac{d}{c} < 1 & \alpha-\beta & 0 \\
\hline
C'_2 & \frac{1}{2} & -2\beta & 0 
& \frac{1}{2}< -\frac{d}{c} < 1 & -2\beta & 0 \\
\hline
D'_2 & \frac{1}{2} & 0 & \alpha-\beta 
& \frac{1}{2}< -\frac{d}{c} < 1 & \alpha-\beta & \alpha-\beta \\
\hline
A'_3 & - & \beta & \beta & - & \beta & \beta \\
\hline
B'_3 & - & -(\alpha-\beta) & 0 & - & -(\alpha-\beta) & 0 \\
\hline
\end{array}$\\
Table 4.1\\
\end{center}

 When $t=1$,
a minimal edge-path contains an edge labeled $C$
only in the $(w+1)$th quadrilateral of the first sequence of $Q_{p/q}$
(the first quadrilateral of the second sequence).
 In fact, such an edge-path is only $c_2$ and $d_2$.
 The transformation
$\left( \begin{array}{cc}
1 & 0 \\
r-1 & 1 \\
\end{array} \right)$
brings the edge $C_0$
to the edge $C'$ oriented from $1/r$ to $0/1$.
 Table 4.2 is derived from Table 2.2,
and shows contribution of $C'$
to the intersection numbers $i_1$ and $i_2$.
 Note that $r \ge 3$.

\begin{center}
$\begin{array}{|c|c|c|c|}
\hline
 & -\frac{d}{c} & i_1 & i_2 \\
\hline
C' & - & 2(\beta-n) & 2n \\
\hline
\end{array}$\\
Table 4.2\\
\end{center}


 For the $j$-th quadrilateral
in the last sequence of quadrilaterals $Q_{p/q}$ with $s>0$,

$E'_0 = g_0 (E_0)$,
$g_0 =
\left( \begin{array}{cc}
-(2j-3) & 2j-2 \\
-\{ (2j-3)r+1 \} & (2j-2)r+1 \\
\end{array} \right)
\in G$,
$-\dfrac{d}{c} = \dfrac{(2j-2)r+1}{(2j-3)r+1}$,

$E'_1 = g_1 (E_0)$,
$g_1 =
\left( \begin{array}{cc}
-(2j-3) & 1 \\
-\{ (2j-3)r+1 \} & r \\
\end{array} \right)
\in G$,
$-\dfrac{d}{c} = \dfrac{r}{(2j-3)r+1}$,

$E'_2 = g_2 (E_0)$,
$g_2 =
\left( \begin{array}{cc}
2j-1 & -1 \\
(2j-1)r+1 & -r \\
\end{array} \right)
\in G$,
$-\dfrac{d}{c} = \dfrac{r}{(2j-1)r+1}$, and

$E'_3 = g_3 (E_0)$,
$g_3 =
\left( \begin{array}{cc}
2j-1 & 2j-2 \\
(2j-1)r+1 & (2j-2)r+1 \\
\end{array} \right)
\in G$,
$-\dfrac{d}{c} = - \dfrac{(2j-2)r+1}{(2j-1)r+1}$

 Hence, when $1 \le t \le \infty$ ($1 < t \le \infty$ for $C'_k$),
we obtain Table 4.3 of contribution
to the intersection numbers $i_1$ and $i_2$
for the edges of the $j$-th quadrilateral
in the second sequence of $Q_{p/q}$ with $s>0$.
 Note that the numbers of the table are independent of $r$.

\begin{center}
$\begin{array}{|c|c|c|c|c|c|c|}
\hline
 & \multicolumn{3}{c|}{j=1} & \multicolumn{3}{c|}{j \ge 2} \\
\hline
 & -\frac{d}{c} & i_1 & i_2 & -\frac{d}{c} & i_1 & i_2 \\
\hline
A'_0 & - & \beta & \beta 
& 1<-\frac{d}{c}<\infty & -\beta & -\beta \\
\hline
B'_0 & - & -(\alpha-\beta) & 0 
& 1<-\frac{d}{c}<\infty & \alpha-\beta & 0 \\
\hline
C'_0 & - & 0 & 2\beta 
& 1<-\frac{d}{c}<\infty & 0 & 2\beta \\
\hline
D'_0 & - & -(\alpha-\beta) & \alpha-\beta 
& 1<-\frac{d}{c}<\infty & \alpha-\beta & \alpha-\beta \\
\hline
A'_1 & - & \beta & \beta 
& 0<-\frac{d}{c}<1 & -\beta & -\beta \\
\hline
B'_1 & - & -(\alpha-\beta) & 0 
&  0<-\frac{d}{c}<1 & \alpha-\beta & 0 \\
\hline
A'_2 & \frac{1}{2}<-\frac{d}{c}<1 & -\beta & -\beta 
& 0<-\frac{d}{c}<\frac{1}{2} & -\beta & -\beta \\
\hline
B'_2 & \frac{1}{2}<-\frac{d}{c}<1 & \alpha-\beta & 0
& 0<-\frac{d}{c}<\frac{1}{2} & \alpha-\beta & 0 \\
\hline
C'_2 & \frac{1}{2}<-\frac{d}{c}<1 & -2\beta & 0 
& 0<-\frac{d}{c}<\frac{1}{2} & -2\beta & 0 \\
\hline
D'_2 & \frac{1}{2}<-\frac{d}{c}<1 & \alpha-\beta & \alpha-\beta 
& 0<-\frac{d}{c}<\frac{1}{2} & -(\alpha-\beta) & \alpha-\beta \\
\hline
A'_3 & - & \beta & \beta & - & \beta & \beta \\
\hline
B'_3 & - & -(\alpha-\beta) & 0 & - & -(\alpha-\beta) & 0 \\
\hline
\end{array}$\\
Table 4.3\\
\end{center}


 For the $j$-th quadrilateral
in the last sequence of quadrilaterals $Q_{p/q}$ with $s<0$,

$E'_0 = h_0 (E_0)$,
$h_0 =
\left( \begin{array}{cc}
-(-2j+1) & -2j+2 \\
-\{ (-2j+1)r+1 \} & (-2j+2)r+1 \\
\end{array} \right)
\in G$,
$-\dfrac{d}{c} = \dfrac{(-2j+2)r+1}{(-2j+1)r+1}$,

$E'_1 = h_1 (E_0)$,
$h_1 =
\left( \begin{array}{cc}
-(-2j+1) & 1 \\
-\{ (-2j+1)r+1 \} & r \\
\end{array} \right)
\in G$,
$-\dfrac{d}{c} = \dfrac{r}{(-2j+1)r+1}$,

$E'_2 = h_2 (E_0)$,
$h_2 =
\left( \begin{array}{cc}
-2j+3 & -1 \\
(-2j+3)r+1 & -r \\
\end{array} \right)
\in G$,
$-\dfrac{d}{c} = \dfrac{r}{(-2j+3)r+1}$, and

$E'_3 = h_3 (E_0)$,
$h_3 =
\left( \begin{array}{cc}
-2j+3 & -2j+2 \\
(-2j+3)r+1 & (-2j+2)r+1 \\
\end{array} \right)
\in G$,
$-\dfrac{d}{c} = - \dfrac{(-2j+2)r+1}{(-2j+3)r+1}$

 Hence, when $1 \le t \le \infty$ ($1 < t \le \infty$ for $C'_k$),
we obtain Table 4.4 of contribution
to the intersection numbers $i_1$ and $i_2$
for the edges of the $j$-th quadrilateral
in the second sequence of $Q_{p/q}$ with $s<0$.
 Note that the numbers of the table are independent of $r$.

\begin{center}
$\begin{array}{|c|c|c|c|c|c|c|}
\hline
 & \multicolumn{3}{c|}{j=1} & \multicolumn{3}{c|}{j \ge 2} \\
\hline
 & -\frac{d}{c} & i_1 & i_2 & -\frac{d}{c} & i_1 & i_2 \\
\hline
A'_0 & - & \beta & \beta 
& \frac{1}{2}<-\frac{d}{c}<1 & -\beta & -\beta \\
\hline
B'_0 & - & -(\alpha-\beta) & 0 
& \frac{1}{2}<-\frac{d}{c}<1 & \alpha-\beta & 0 \\
\hline
C'_0 & - & 0 & 2\beta 
& \frac{1}{2}<-\frac{d}{c}<1 & -2\beta & 0 \\
\hline
D'_0 & - & -(\alpha-\beta) & \alpha-\beta 
& \frac{1}{2}<-\frac{d}{c}<1 & \alpha-\beta & \alpha-\beta \\
\hline
A'_1 & - & \beta & \beta 
& - & \beta & \beta \\
\hline
B'_1 & - & -(\alpha-\beta) & 0 
&  - & -(\alpha-\beta) & 0 \\
\hline
A'_2 & \frac{1}{2}<-\frac{d}{c}<1 & -\beta & -\beta 
& - & \beta & \beta \\
\hline
B'_2 & \frac{1}{2}<-\frac{d}{c}<1 & \alpha-\beta & 0
& - & -(\alpha-\beta) & 0 \\
\hline
C'_2 & \frac{1}{2}<-\frac{d}{c}<1 & -2\beta & 0 
& - & 0 & 2\beta \\
\hline
D'_2 & \frac{1}{2}<-\frac{d}{c}<1 & \alpha-\beta & \alpha-\beta 
& - & -(\alpha-\beta) & \alpha-\beta \\
\hline
A'_3 & 0<-\frac{d}{c}\le\frac{1}{2} & -\beta & -\beta 
& - & \beta & \beta \\
\hline
B'_3 & 0<-\frac{d}{c}\le\frac{1}{2} & \alpha-\beta & 0 
& - & -(\alpha-\beta) & 0 \\
\hline
\end{array}$\\
Table 4.4\\
\end{center}

\section{Euler characteristic}

 We give a formula for the Euler characteristic $\chi(F)$ in this section.
 We consider only the case of $1 \le t \le \infty$.
 (In the case of $0 \le t \le 1$,
we simply exchange $\alpha$ and $\beta$.)

 Let $e_1, e_2, \cdots, e_k$ be the edges of the minimal edge-path.
 Let $F_i$ be the part of $F$ carried by $\Sigma_{e_i}$.
 Then $\chi(F) = (\Sigma_{i=1}^k \chi(F_i)) - (k-1)(\alpha+\beta)$
since $F_i$ and $F_{i+1}$ are glued together
along $\alpha + \beta$ arcs.
 Each $\chi(F_i)$ is as in Table 5.1
according to the label of the edge $e_i$.

\vspace{2mm}
\begin{center}
$\begin{array}{|c|c|}
\hline
{\rm label} & \chi(\Sigma_{e_i}) \\
\hline
A & \alpha \\
\hline
B & 2\beta + \frac{\alpha - \beta}{2} \\
\hline
C & \alpha \\
\hline
D & 2 \beta \\
\hline
\end{array}$\\
Table 5.1\\
\end{center}

\begin{figure}[htbp]
\begin{center}
\includegraphics[width=7cm]{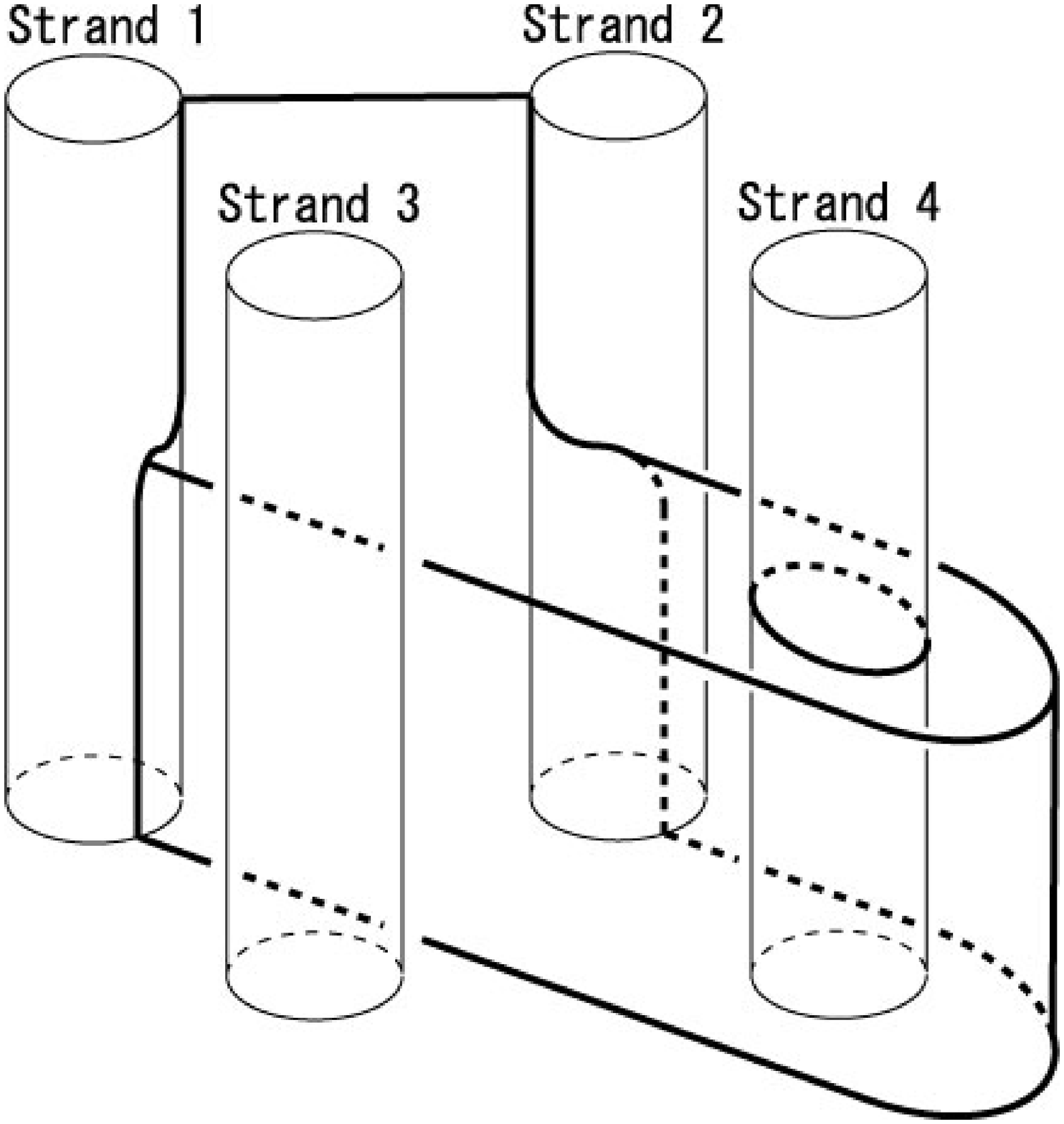}
\end{center}
\caption{}
\label{fig:SigmaD3}
\end{figure}

\begin{figure}[htbp]
\begin{center}
\includegraphics[width=7cm]{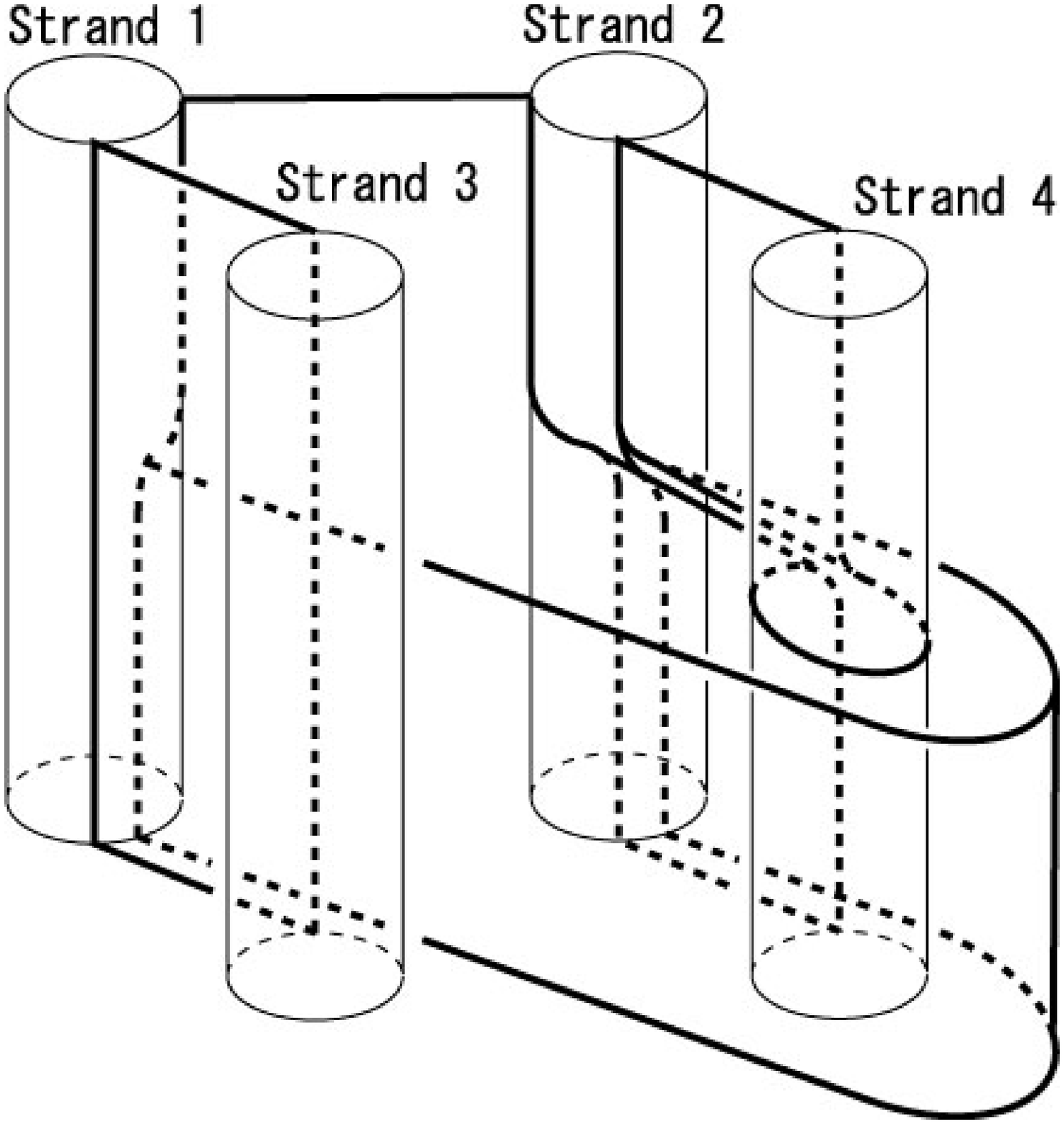}
\end{center}
\caption{}
\label{fig:SigmaD4}
\end{figure}

 In fact, for the label $A$,
$F_i$ is a disjoint union of $(\alpha-\beta) + \beta = \alpha$ discs.
 For the label $B$,
$F_i$ is a disjoint union of
$2\beta + \frac{\alpha - \beta}{2}$ discs.
 For the label $C$,
$F_i$ is a disjoint union of $(\alpha-\beta) + \beta = \alpha$ discs.
 For the label $D$,
$F_i$ is a disjoint union of $2 \beta$ discs
and $0$ or more annuli (discs punctured by $L_2$).
 Precisely, $F_i$ has
$\alpha$ annuli when $\beta = 0$,
and no annuli when $\beta \ne 0$.
 See Figures \ref{fig:SigmaD3} and \ref{fig:SigmaD4},
where $F_i$ with label $D$
with $\alpha = 1, \beta = 0, n = 1$
and $\alpha = 2, \beta = 1, n = 1$ respectively
are described.

\section{Generalized genus}

 Let $F$ be a compact $2$-manifold.
 $F$ may be orientable or non-orientable,
and connected or disconnected.
 $g'(F)$ denotes the {\it generalized genus}
defined by $\chi(F) = 2-2g'(F)-b(F)$,
where $\chi(F)$ is the Euler characteristic of $F$,
and $b(F)$ the number of boundary circles.

\begin{lemma}\label{lem:genus}
 If $g'(F)=0$,
then $F$ has a planar surface component
or a projective plane with holes component.

 If $g'(F)=1$,
then $F$ has either
a planar surface component,
a projective plane with holes component,
a torus with holes component
or a Klein bottle with holes component.
\end{lemma}

\begin{proof}
 For a connected orientable surface $F_0$,
the usual genus is calculated by
$g(F_0) = \{ 2 - \chi(F_0) - b(F_0) \} / 2$.
 Hence $b(F_0) \le -\chi(F_0)$ when $g(F_0) \ge 1$,
and $b(F_0) \le -\chi(F_0) -2$ when $g(F_0) \ge 2$.

 For a connected non-orientable surface $F_0$,
set $t(F_0) = 2 - \chi(F_0) -b(F_0)$.
 Then $t(F_0) \ge 1$.
 $F_0$ is a projective plane (with holes) if $t(F_0)=1$.
 $F_0$ is a Klein bottle (with holes) if $t(F_0)=2$.
 Moreover, $b(F_0) \le -\chi(F_0)$ when $t(F_0) \ge 2$,
and $b(F_0) \le -\chi(F_0) -1$ when $t(F_0) \ge 3$.


 The calculation below shows the latter half of this lemma.
 (The former half follows similarly.)
 Set $F = F_1 \cup ... \cup F_n$,
where $F_i$ is a connected component for $1 \le i \le n$.
 If $g(F_i) \ge 2$ for each orientable component $F_i$ of $F$
and if $t(F_j) \ge 3$ for each non-orientable component $F_j$,
then $b(F) = \Sigma_k b(F_k) \le \Sigma_k (-\chi(F_k)-1) = -\chi(F)-n$.
 Then $g'(F) = \{ 2 -\chi(F)-b(F) \} /2 \ge \{ 2+ (b(F)+n)-b(F)\} /2
= (n+2)/2 = n/2 + 1 > 1$.
\end{proof}

\begin{remark}
 Let $F$ be a compact $2$-manifold
properly embedded in an orientable $3$-manifold $M$.
 Let $F_1, \cdots, F_n$ be orientable connected components of $F$,
and $P_1, \cdots, P_m$ non-orientable ones.
 A regular neighbourhood $N(P_i)$ of each $P_i$
is a twisted $I$-bundle over $P_i$.
 Let $\tilde{P_i}$ be the frontier surface of $N(P_i)$,
that is, $\tilde{P_i} = {\rm cl}\,(\partial N(P_i) - \partial M)$,
which is a connected orientable surface,
and called the {\it double} of $P_i$ in this paper.
 Then, an easy calculation shows
$g'(F) = 1 + \Sigma_{i=1}^n (g(F_i)-1)
 + \frac{1}{2} \Sigma_{j=1}^m (g(\tilde{P_i})-1)$.

 Note
that the double of a projective plane with holes is a planar surface,
and that the double of a Klein bottle with holes is a torus with holes.
\end{remark}


\section{Surface $1$-$6$ with $s>0$ and $1 \le t \le \infty$}

 In this section,
we calculate boundary slope of the surface $F$
corresponding to the edge-path $c_{16}$ with $s>0$ and $1 \le t \le \infty$.
 We are going to obtain the slope of the preferred longitude
by substituting $1$ for $\alpha$, and $0$ for $\beta$.

 Table 7.1 shows the calculation which uses Tables 4.1 and 4.3.
 We can find from the left side
the columns of vertices, quadrilaterals,
labels of edges, $i_1$, $i_2$ and $\chi(F_i)$.
 In the second column,
the sign of $-A'_2$ means
that the orientations of $e_i$ and $g(A_0)$ don't match.
 In the $5$--$7$th rows, $2 \le i \le w$.
 In the third row from the bottom, $2 \le j \le u+1$.

 Thus, from Table 7.1, the $\lq\lq$slope" on $\partial N(L_1)$ is
$(\alpha, (w-u)(\alpha-\beta) + (2w+1)\beta)$
and that on $\partial N(L_2)$ is 
$(\beta, (w+u+1)(\alpha-\beta) + (2w+1)\beta)$,
where the first coordinate is the longitudinal entry,
and the second coordinate is
the $\lq\lq$meridional" entry with respect to the unusual longitude $\lambda_1$.
 To obtain the real slope,
we need to know the $\lq\lq$slope" of the preferred longitude,
and divide the entries by their greatest common measure.
 $\partial F \cap \partial N(L_1)$ has
$GCM(\alpha, (w-u)(\alpha-\beta) + (2w+1)\beta)
 = GCM(\alpha, (w-u)\alpha + (w+u+1)\beta)
 = GCM(\alpha, (w+u+1)\beta)$ circles,
and $\partial F \cap \partial N(L_2)$ has
$GCM(\beta, (w+u+1)(\alpha-\beta) + (2w+1)\beta)
= GCM(\beta, (w+u+1)\alpha + (w-u)\beta)
= GCM(\beta, (w+u+1)(\alpha-\beta))$ circles.

 Substituting $1$ for $\alpha$ and $0$ for $\beta$,
we obtain the $\lq\lq$slope" of the preferred longitude of $L_1$,
which is $(1, w-u)$.
 In this case,
the minimal edge-path is that of $c_6$ in $D_{\infty}$,
and is composed only of $w+u+1$ edges labeled $D$.
 The surface carried by $\Sigma_D$ with $\alpha=1$ and $\beta=0$
is a disc punctured by $L_2$ once.
 See Figure \ref{fig:SigmaD3}.
 Hence $F$ is a disc bounded by $L_1$ punctured by $L_2$ $w+u+1$ times.

 Thus the slopes with respect to the preferred longitude can be obtained from
$(\alpha, (w-u)\alpha + (w+u+1)\beta -(w-u)\alpha)
 = (\alpha, (w+u+1)\beta)$
on $\partial N(L_1)$, and 
$(\beta, (w+u+1)\alpha + (w-u)\beta -(w-u)\beta)
 = (\beta, (w+u+1)\alpha)$
on $\partial N(L_2)$
by dividing by the G.C.M. of the longitudinal entry and the meridional entry.

\vspace{2mm}
\begin{center}
$\begin{array}{|c|c|c|c|c|c|}
\hline
{\rm vertex} & {\rm quadrilateral} & {\rm label} & i_1 & i_2 & \chi(F_i) \\
\hline
\frac{1}{0}\sim & i=1 & A'_1 
 & 0        & 0 & \alpha \\
\ & & D'_2 
 & 0        & \alpha-\beta & 2 \beta \\
\sim\frac{1}{2} & & -A'_2 
 & \beta        & \beta & \alpha \\
\hline
\frac{1}{2i-2}\sim & 2 \le i \le w & A'_1 
 & \beta        & \beta & \alpha \\
\ & & D'_2 
 & \alpha-\beta  & \alpha-\beta & 2\beta \\
\sim\frac{1}{2i} & & -A'_2 
 & \beta        & \beta & \alpha \\
\hline
\frac{1}{r-1}\sim & i = w+1 & A'_1 
 & \beta        & \beta & \alpha \\
\ & & D'_2 
 & \alpha-\beta  & \alpha-\beta & 2\beta \\
\hline
\ & 2 \le j \le u+1 & D'_2 
 & -(\alpha-\beta)  & \alpha-\beta & 2\beta \\
\hline
\sim [r,s] & j= u+1 & -A'_2 
 & \beta        & \beta & \alpha \\
\hline
\ & {\rm total} & \alpha-\beta \ {\rm entry}
 & w-u & w+u+1 & 2(w+1) \\
\cline{3-6}
\ & & \beta \ {\rm entry}
 & 2w+1 & 2w+1 & 2(2w+u+2) \\
\hline
\end{array}$ \\
Table 7.1
\end{center}

 We can calculate Euler characteristic by
\newline
$\chi(F) = 2(w+1)(\alpha-\beta)+2(2w+u+2)\beta -\{ (3w+u+3)-1\} (\alpha+\beta)
= -(w+u)(\alpha-\beta) - 2w\beta$.

 Thus the generalized genus is
\newline
$g'(F) = 
\{ (w+u)(\alpha-\beta) + 2w\beta
-GCM(\alpha, (w+u+1)\beta)-GCM(\beta, (w+u+1)(\alpha-\beta))
+2 \} /2$.

\section{Surface $2$-$6$ with $s<0$ and $1<t \le \infty$}

 In this section,
we calculate boundary slope of the surface $F$
corresponding to the edge-path $d_{26}$ with $s<0$ and $1<t \le \infty$.
 We set $s = 2u+1 = -(2u'+1)$.
 We are going to obtain the slope of the preferred longitude
by substituting $1$ for $\alpha$, and $0$ for $\beta$.

 Table 8.1 shows the calculation which uses Tables 4.1 and 4.4.
 In the $3$rd row, $2 \le i \le w$.
 In the $5$th row, $2 \le j \le u'+1$.

 Thus the $\lq\lq$slope" on $\partial N(L_1)$ is
$(\alpha, (w+u'-1)(\alpha - \beta) - \beta)$
and that on $\partial N(L_2)$ is
$(\beta, -(w+u')(\alpha - \beta) -3\beta)$.
 $\partial F \cap \partial N(L_1)$ has
$b_1 = GCM(\alpha, (w+u'-1)(\alpha - \beta) - \beta)
= GCM(\alpha, (w+u')\beta)$ circles,
and $\partial F \cap \partial N(L_2)$ has
$b_2 = GCM(\beta, -(w+u')(\alpha - \beta) -3 \beta)
= GCM(\beta, (w+u')(\alpha -\beta))$ circles.

 Substituting $1$ for $\alpha$ and $0$ for $\beta$,
we obtain the $\lq\lq$slope" of the preferred longitude of $L_1$,
which is $(1, w+u'-1)$.
 In this case,
$F$ is a disc bounded by $L_1$ punctured by $L_2$ $w+u'$ times.

 Thus the slopes with respect to the preferred longitude can be obtained from
$(\alpha, (w+u'-1)(\alpha-\beta) -\beta -(w+u'-1)\alpha)
= (\alpha, -(w+u')\beta)$
on $\partial N(L_1)$, and 
$(\beta, -(w+u')(\alpha-\beta) -3\beta -(w+u'-1)\beta) 
= (\beta, -(w+u')\alpha-2\beta)$
on $\partial N(L_2)$
by divided by the G.C.M. of the longitudinal entry and the meridional entry.

\vspace{2mm}
\begin{center}
$\begin{array}{|c|c|c|c|c|c|}
\hline
{\rm vertex} & {\rm quadrilateral} & {\rm label} & i_1 & i_2 & \chi(F_i) \\
\hline
\frac{1}{0}\sim & i=1 & A'_0 
 & 0        & 0 & \alpha \\
\ & & -D'_0 
 & 0        & -(\alpha-\beta) & 2 \beta \\
\hline
\ & 2\le i \le w & -D'_0 
 & \alpha-\beta  & -(\alpha-\beta) & 2\beta \\
\hline
\ & i=w+1 & -C'_0 
 & 0        & -2\beta & \alpha \\
\hline
\ & 2\le j \le u'+1 & -D'_2 
 & \alpha-\beta  & -(\alpha-\beta) & 2\beta \\
\hline
\sim [r,s] & j=u'+1 & -A'_1 
 & -\beta        & -\beta & \alpha \\
\hline
\ & {\rm total} & \alpha-\beta \ {\rm entry}
 &  w+u'-1 & -(w+u') & 3 \\
\cline{3-6}
\ & & \beta \ {\rm entry}
 & -1 & -3 & 2w+2u'+3 \\
\hline
\end{array}$ \\
Table 8.1
\end{center}

 We can calculate Euler characteristic by
\newline
$\chi(F)
= 3(\alpha-\beta)+(2w+2u'+3)\beta -\{ (w+u'+3)-1\} (\alpha+\beta)
= -(w+u'-1)(\alpha - \beta) - \beta$.

 Thus the generalized genus is
\newline
$g'(F) = 
\{ (w+u'-1)(\alpha - \beta) + \beta - b_1 - b_2 + 2 \} /2$,
\newline
where
$b_1 = GCM(\alpha, (w+u')\beta)$ and
$b_2 = GCM(\beta, (w+u')(\alpha - \beta))$.

\section{data of all the essential surfaces for $L([r,s])$}

 We can obtain boundary slopes, Euler characteristic, generalized genus
of the surfaces corresponding to all the minimal paths as below
by similar calculations as in previous two sections.

 The boundary slopes below are
with respect to the ordinary preferred longitudes,
each of which is the boundary of a disc
bounded by a component $L_i$ and punctured by the other component $L_j$.
 But they are not divided
by the G.C.M. of the longitudinal entry and the meridional entry.

 For $L([r,s])$ with $s > 0$,
we obtain Tables 9.1 and 9.2 below,
where $c_2$ is for $t=1$,
$c_{2k}$ is for $1 < t \le \infty$
and the others are for $1 \le t \le \infty$.

 From data of $c_{xy}$ with $x \ne 2$,
we can obtain those for $c_{x}$ by substituting $\beta$ for $\alpha$,
and those for $c_{y}$ by substituting $0$ for $\beta$.
 We omit the data of the surfaces for $0 \le t < 1$
because they are obtained
from those for $1 \le t \le \infty$
by $\pi$-rotation of $L([r,s])$
about the axis $\Gamma(\frac{1}{2}, \frac{1}{2}) \times {\Bbb R}$.
 We omit also the data of the surfaces
corresponding to $c_{38}$, $c_{36}$, $c_{28}$, $c_{27}$
since they are obtained from those for $L([s,r])$
corresponding to $c_{14}$, $c_{16}$, $c_{24}$, $c_{25}$
by reflecting $L([s,r])$ upside down
and transforming every level sphere by
$\left( \begin{array}{cc}
s & -1 \\
sr+1 & -r \\
\end{array} \right)$.


\begin{center}
\begin{tabular}{|c|c|c|c|}
\hline
path & slope on $\partial N(L_1)$ & slope on $\partial N(L_2)$ 
& Euler char. \\
\hline
$c_2$ 
& $(\beta, -(w-u+1)\beta+2n)$ & $(\beta, -(w-u-1)\beta-2n)$ 
& $-\beta$ \\
\hline
$c_{14}$ 
& $(\alpha, (u+1)\alpha + w\beta)$ & $(\beta, w\alpha + (u+1)\beta)$ 
& $-w(\alpha+\beta)$ \\
\hline
$c_{16}$ 
& $(\alpha, (w+u+1)\beta)$ & $(\beta, (w+u+1)\alpha)$ 
& $-(w+u)(\alpha-\beta)-2w\beta$ \\
\hline
$c_{24}$ 
& $(\alpha, (u+1)\alpha-w\beta)$ & $(\beta, -w\alpha+(u-1)\beta)$ 
& $-w(\alpha-\beta)-\beta$ \\
\hline
$c_{25}$ 
& $(\alpha, -(w-u-1)\alpha)$ & $(\beta, -(w-u+1)\beta)$ 
& $-\alpha$ \\
\hline
\end{tabular}\\
Table 9.1
\end{center}

 Let $b_i$ be the number of boundary circles on $\partial N(L_i)$
for $i=1$ and $2$.
 Set $b = b_1 + b_2$.

\begin{center}
\begin{tabular}{|c|c|}
\hline
path & generalized genus \\
\hline
$c_2$ &
$g' = (\beta - b + 2)/2$,
$\ \ b_1 = b_2 = GCM(\beta, 2n)$ \\
\hline
$c_{14}$ &
$g' = \{ w(\alpha+\beta)-b+2 \} /2$,
$\ \ b_1 = GCM(\alpha, w\beta)$, $\ b_2 = GCM(\beta, w\alpha)$ \\
\hline
$c_{16}$ &
$\begin{array}{c}
g' = \{ (w+u)(\alpha-\beta)+2w\beta-b+2 \} /2, \\
b_1 = GCM(\alpha, (w+u+1)\beta), \ b_2 = GCM(\beta, (w+u+1)\alpha) \\
\end{array}$ \\
\hline
$c_{24}$ &
$g' = \{ w(\alpha-\beta)+\beta -b +2 \} /2$,
$\ \ b_1 = GCM(\alpha, w\beta), \ b_2 = GCM(\beta, w\alpha)$ \\
\hline
$c_{25}$ &
$g' = (-\beta+2)/2$,
$\ \ $ ($b_1 = \alpha,\ b_2 = \beta$) \\
\hline
\end{tabular}\\
Table 9.2
\end{center}

 For $L([r,s])$ with $s < 0$,
we obtain Tables 9.3 and 9.4 below,
where $d_2$ is for $t=1$,
$d_{2k}$ is for $1 < t \le \infty$
and the others are for $1 \le t \le \infty$.

 We omit also the data of the surfaces
corresponding to $d_{38}$, $d_{36}$, $d_{28}$
since they are obtained from those for $L([-s,-r])$
corresponding to $d_{14}$, $d_{16}$, $d_{24}$
by $\pi$-rotating $L([-s,-r])$ about a horizontal line parallel to the $x$-axis
and transforming every level sphere by
$\left( \begin{array}{cc}
s & -1 \\
sr+1 & -r \\
\end{array} \right)$.

\begin{center}
\begin{tabular}{|c|c|c|c|}
\hline
path & slope on $\partial N(L_1)$ & slope on $\partial N(L_2)$ 
& Euler char. \\
\hline
$d_2$ 
& $(\beta, -(w-u+1)\beta+2n)$ & $(\beta, -(w-u-1)\beta-2n)$ 
& $-\beta$ \\
\hline
$d_{06}$ 
& $(\alpha, (w+u+1)\beta)$ & $(\beta, (w-u-1)\alpha)$
& $-(w-u-2)(\alpha+\beta)$ \\
\hline
$d_{14}$ 
& $(\alpha, (u+1)\alpha + w\beta)$ & $(\beta, w\alpha + (u+1)\beta)$ 
& $-w(\alpha+\beta)$ \\
\hline
$d_{16}$ 
& $(\alpha, (w+u+1)\beta)$ & $(\beta, (w+u+1)\alpha)$ 
& $-(w-u-2)(\alpha-\beta)-2w\beta$ \\
\hline
$d_{24}$ 
& $(\alpha, (u+1)\alpha-w\beta)$ & $(\beta, -w\alpha+(u-1)\beta)$ 
& $-w(\alpha-\beta)-\beta$ \\
\hline
$d_{25}$ 
& $(\alpha, -(w-u-1)\alpha)$ & $(\beta, -(w-u+1)\beta)$ 
& $-\alpha$ \\
\hline
$d_{26}$ 
& $(\alpha, -(w-u-1)\beta)$ & $(\beta, -(w-u-1)\alpha-2\beta)$ 
& $-(w-u-2)(\alpha-\beta) + \beta$ \\
\hline
$d_{27}$ 
& $(\alpha, -(w-u+1)\alpha)$ & $(\beta, -(w-u-1)\beta)$ 
& $-\alpha$ \\
\hline
\end{tabular}\\
Table 9.3
\end{center}

\begin{center}
\begin{tabular}{|c|c|}
\hline
path & generalized genus \\
\hline
$d_2$ &
$g' = (\beta - b + 2)/2$,
$\ \ b_1 = b_2 = GCM(\beta, 2n)$ \\
\hline
$d_{14}$ &
$g' = \{ w(\alpha+\beta) -b +2 \} /2$,
$\ \ b_1 = GCM(\alpha, w\beta), \ b_2 = GCM(\beta, w\alpha)$ \\
\hline
$d_{06}$ &
$\begin{array}{c}
g' = \{ (w-u-2)(\alpha+\beta) -b +2 \} /2, \\
b_1 = GCM(\alpha, (w+u+1)\beta), \ b_2 = GCM(\beta, (w-u-1)\alpha) \\
\end{array}$ \\
\hline
$d_{16}$ &
$\begin{array}{c}
g' = \{ (w-u-2)(\alpha-\beta)+2w\beta -b +2 \} /2, \\
b_1 = GCM(\alpha, (w+u+1)\beta), \ b_2 = GCM(\beta, (w+u+1)\alpha) \\
\end{array}$ \\
\hline
$d_{24}$ &
$g' = \{ w(\alpha-\beta)+\beta -b +2 \} /2$,
$\ \ b_1 = GCM(\alpha, w\beta), \ b_2 = GCM(\beta, w\alpha)$ \\
\hline
$d_{25}$ &
$g' = (-\beta+2)/2$,
$\ \ $ ($b_1=\alpha, \ b_2 = \beta$)  \\
\hline
$d_{26}$ &
$\begin{array}{c}
g' = \{ (w-u-2)(\alpha-\beta) + \beta -b + 2 \} /2, \\
b_1 = GCM(\alpha, (w-u-1)\beta),\ b_2 = GCM(\beta, (w-u-1)\alpha) \\
\end{array}$ \\
\hline
$d_{27}$ &
$g' = (-\beta+2)/2$,
$\ \ $ ($b_1=\alpha, \ b_2 = \beta$)  \\
\hline
\end{tabular}\\
Table 9.4
\end{center}

\section{Proof of necessity}

 We prove that, if $L([r,s])[\gamma_1, \gamma_2]$ is reducible,
then the slopes $\gamma_1, \gamma_2$ are as in Theorem \ref{thm:red}.

 If a result of non-trivial Dehn surgery $L([r,s])[\gamma_1, \gamma_2]$
is reducible,
then it is well-known
that the exterior of $L([r,s])$ contains
an incompressible and boundary incompressible planar surface $P$
with non empty boundary circles of slope $\gamma_i$ on $\partial N(L_i)$
for $i=1$ and $2$.
 Because $\gamma_1, \gamma_2$ are non-meridional,
the surface $P$ is meridionally incompressible.
 Hence $P$ is carried by a branched surface $\Sigma_{\gamma}$
corresponding to a minimal edge-path $\gamma$ in $Q_{[r,s]}$
as shown in \cite{FH}.
 So, we should know
when generalized genus listed in the previous section
can be zero.
 In fact,
it can be zero only for $c_2, d_2, c_{25}, d_{25}, d_{26}$ and $d_{27}$
by calculation as below.

{\bf Calculation for $c_2$ and $d_2$}.
\newline
 The generalized genus is zero only when $\beta=2$ and $n=0,1$ or $2$
by the calculation below.
 Moreover, by Tables 9.1 and 9.3,
we can calculate numbers of boundary circles and slopes.
 When $n=2$,
$2$ circles of slope $(1, -w+u+1)$ on $\partial N(L_1)$ and 
$2$ circles of slope $(1, -w+u-1)$ on $\partial N(L_2)$,
which is (1)(a) in Theorem \ref{thm:red}.
 When $n=1$,
$2$ circles of slope $(1, -w+u)$
both on $\partial N(L_1)$ and on $\partial N(L_2)$,
which is (1)(b) in Theorem \ref{thm:red}.
 When $n=0$, 
$2$ circles of slope $(1, -w+u-1)$ on $\partial N(L_1)$ and 
$2$ circles of slope $(1, -w+u+1)$ on $\partial N(L_2)$,
which is (1)(a) in Theorem \ref{thm:red}.

 $g'=\{ \beta -2GCM(\beta, 2n) +2 \} /2$ as in the lists.
 If $g'=0$, then $\beta = 2GCM(\beta, 2n) - 2 = 2(GCM(\beta, 2n) - 1)$.
 Hence $\beta = 2 \beta'$, for some positive integer $\beta'$.
 So, we have $\beta' = 2 GCM(\beta', n)-1 \cdots$ (i).
 Recall that $0 \le n \le \beta = 2\beta'$.
 If $GCM(\beta', n) \le \beta'/2$
then $\beta' = 2 GCM(\beta', n)-1 \le 2 \frac{\beta'}{2}-1 = \beta'-1$,
which is a contradiction.
 Hence $GCM(\beta', n) = \beta'$,
and $\beta'=1$ by (i).

{\bf Calculation for $c_{25}$ and $d_{25}$}.
\newline
 The generalized genus $(-\beta+2)/2$ is zero only when $\beta=2$.
 Moreover, by Tables 9.1 and 9.3, the surface with $\beta=2$ has
$\alpha$ boundary circles of slope $(1, -w+u+1)$ on $\partial N(L_1)$ and 
$2$ boundary circles of slope $(1, -w+u-1)$ on $\partial N(L_2)$,
which is (1)(a) in Theorem \ref{thm:red}.

{\bf Calculation for $d_{27}$}.
\newline
 The generalized genus $(-\beta+2)/2$ is zero only when $\beta=2$.
 Moreover, by Table 9.3, the surface with $\beta=2$ has
$\alpha$ boundary circles of slope $(1, -w+u-1)$ on $\partial N(L_1)$ and 
$2$ boundary circles of slope $(1, -w+u+1)$ on $\partial N(L_2)$,
which is (1)(a) in Theorem \ref{thm:red}.

\begin{lemma}\label{lem:AlphaBeta}
 Let $F$ be a surface properly embedded in the exterior of $L([r,s])$
such that the union of boundary circles $\partial F$ goes around longitudinally
$\alpha$ times on $\partial N(L_1)$
and $\beta$ times on $\partial N(L_2)$.
 Suppose $\alpha \ge \beta > 0$.
 If the generalized genus $g'$ of $F$ is zero,
then $\chi \ge -(\alpha + \beta) + 2$,
where $\chi$ is the Euler characteristics of $F$.
\end{lemma}

\begin{proof}
 Let $b$ be the number of boundary circles of $F$.
 We have $b \le \alpha + \beta$.
 Recall $g' = \{ -\chi -b +2 \} /2$.
 If $g' = 0$,
then $\chi = -b + 2 \ge - (\alpha + \beta) + 2$.
\end{proof}

%

{\bf Calculations for $c_{14}, c_{16}, d_{14}, d_{06}, d_{16}$}.
\newline
 When $1 \le t < \infty$,
Lemma \ref{lem:AlphaBeta} shows
that generalized genera of the surfaces
corresponding to $c_{14}, c_{16}, d_{14}, d_{06}, d_{16}$
are never zero.
 For example,
the Euler characteristics of the surface corresponding to $d_{16}$ is
$-(w-u-2)(\alpha-\beta)-2w\beta 
\le -1 \cdot (\alpha-\beta) - 2 \cdot 1 \cdot \beta < -(\alpha+\beta) + 2$.
 Note that $u \le -2$ for $d_{kl}$ since $s \le -3$.

 When $t = \infty$, we have $\beta = 0$.
 The surface corresponding to $c_4$ or $d_4$ has
$w \alpha$ boundary circles of meridional slope on $\partial N(L_2)$,
and the surface corresponding to $c_6$ or $d_6$ has
$(w+u+1) \alpha$ boundary circles of meridional slope on $\partial N(L_2)$.
 Since we are considering non-trivial Dehn surgery on $L([r,s])$,
we can skip this case.

{\bf Calculation for $d_{26}$}.
\newline
 Generalized genus of the surface corresponding to $d_{26}$ is zero
only when $w=1, u=-2, \alpha=4, \beta=2$.
 In this case, the surface has
$4$ boundary circles of slope $(1, -1)$ on $\partial N(L_1)$
and $2$ boundary circles of slope $(1, -6)$ on $\partial N(L_2)$,
which is (2) in Theorem \ref{thm:red}.
 The calculation is a little harder.
 So, we give it in Appendix A.

{\bf Calculation for $c_{24}$, $d_{24}$}.
\newline
 Generalized genera of the surfaces corresponding to $c_{24}, d_{24}$
are never zero.
 The calculation is similar to and much easier than that for $d_{26}$.
 So, we omit it.

\section{Proof of sufficiency}

 In this section,
we prove Theorems \ref{thm:torus}, \ref{thm:cable} in section 1
and Theorem \ref{thm:TorusCable} below,
which shows sufficiency of Theorem \ref{thm:red}.

 Let $K$ be a knot in a $3$-manifold $M$ .
 $K$ is called a {\it trivial knot} if it bounds a disc in $M$.
 $K$ is called a {\it core knot}
if its exterior cl\,$(M-N(K))$ is a solid torus.
 $K$ is called a {\it torus knot}
if it is entirely contained in a Heegaard splitting torus of $M$.
 $K$ is called a {\it cable knot}
if it is isotoped in $\partial N(K')$ for some non-trivial non-core knot $K'$
so that $K$ winds $2$ or more times longitudinally
in the solid torus $N(K')$.

 If we perform Dehn surgery of slope $\gamma$
on one component of a $2$-bridge link $L([r,s])$,
then the other component forms a knot
in a lens space, $S^1 \times S^2$ or $S^3$.
 We let $L([r,s])[\gamma]$ denote it.
 We assume that the surgery is performed on $L_2$ rather than $L_1$,
following the proof of Theorem 5.1 in \cite{W} by Wu.

 Let $L, M$ be the preferred longitude and the meridian of $N(L_2)$
before the Dehn surgery
with the slope $+1$ being $L+M$.
 Let $C'$ be a core circle of the exterior solid torus $E$ of $L_2$
such that $C'$ is homologous to $M$ in $E$.
 We take $M$ and $L$ as a longitude and meridian system for $N(C')$,
exchanging $L$ and $M$.
 Let $C$ be a core circle of the filled solid torus.
 We take $M$ and $\gamma$ as a longitude and meridian system for $N(C)$.
 Note that $\gamma$ is eventually of integral slope in $N(L_2)$
in the next theorem.

\begin{theorem}\label{thm:TorusCable}
 Let $w, u$ be integers with $w \ge 1$ and either $u \ge 1$ or $u \le -2$.
 Then $L([2w+1, 2u+1])[\gamma]$ is
a non-trivial non-core torus knot or a cable knot
if and only if it is in one of the three cases below
or its mirror image.
\begin{enumerate}
\item
 The knot $L([2w+1, 2u+1])[-w+u]$
is a torus knot in the $(-w+u, 1)$-lens space,
which can be placed in $\partial N(L_2)$
so that it is a $(2w+1, -2)$-cable of $C$,
and a $(2u+1, 2)$-cable of $C'$.
\item
 For some integer $k$,
the knot $L([2w+1, 2u+1])[-w+u +\epsilon]$ with $\epsilon = \pm 1$
is a $(2,k)$-cable of a torus knot $K$
in the $(-w+u+\epsilon, 1)$-lens space,
where $K$ can be placed in $\partial N(L_2)$
so that it is a $(\{(2u+1) + \epsilon \} /2, 1)$-cable of $C'$,
and a $(1, -\{ (2w+1) - \epsilon \} /2 )$-cable of $C$.
 Moreover,
$L([2w+1, 2u+1])[-w+u+\epsilon]$ can be isotoped into $\partial N(L_2)$
to be two parallel copies of $K$
except near a single $-\epsilon$ crossing.
 It is a cable knot,
when $(w, \epsilon) \ne (1,1)$ and $(u, \epsilon) \ne (1,-1), (-2,1)$.
 When $(w, \epsilon) = (1,1)$,
$L([3, 2u+1])[u]$ is a torus knot
which can be placed in $\partial N(L_2)$
so that it is a $(3,-2)$-cable of $C$,
and a $(3u+2, 3)$-cable of $C'$.
 When $(u, \epsilon) = (-2,1)$,
$L([2w+1, -3])[-w-1]$ is a torus knot
which can be placed in $\partial N(L_2)$
so that it is a $(3w+1, -3)$-cable of $C$,
and a $(2, -3)$-cable of $C'$.
\item
 $L([3,-3])[-1]$ is a trefoil knot in $S^3$.
\end{enumerate}
\end{theorem}

 Necessity of Theorem \ref{thm:TorusCable} follows
from that of Theorem \ref{thm:red}.
 Conversely, sufficiency of Theorem \ref{thm:red} follows
from that of Theorem \ref{thm:TorusCable}
since a single Dehn surgery on a torus knot
or a cable knot of a torus knot
yields a reducible manifold.
 See \cite{M} and \cite{G}.

 Arguments similar to the proof of Theorem 5.1 in \cite{W} show sufficiency
of Theorem \ref{thm:TorusCable}, \ref{thm:torus} and \ref{thm:cable}.
 He found the case (2) there.
 The case (3) is well-known.

\begin{figure}[htbp]
\begin{center}
\includegraphics[width=7cm]{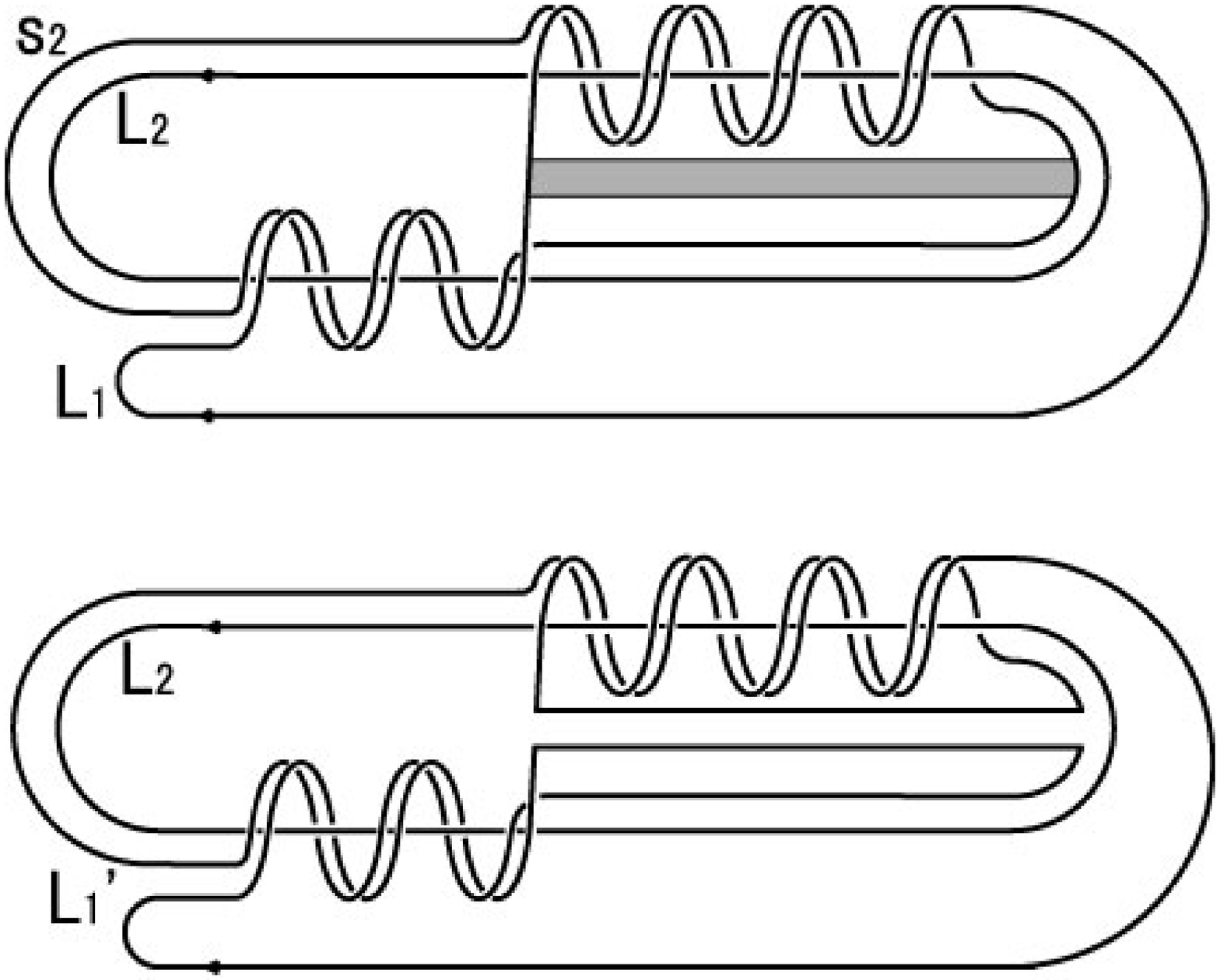}
\end{center}
\caption{}
\label{fig:BandSum2}
\end{figure}

 Let $\cdot$ denote the algebraic intersection number on $\partial N(L_2)$
with $L \cdot M = + 1$ rather than $-1$.

\begin{proof}
 We prove Theorem \ref{thm:TorusCable} (1) and \ref{thm:torus} (1).

 $L([2w+1,2u+1])[\infty,-w+u]$ is the $3$-manifold
obtained from $S^3$ by $(-w+u)$-Dehn surgery on the component $L_2$,
and hence it is the $(-w+u, 1)$-lens space.
 It is again the $3$-sphere when $-w+u = \pm 1$, 
and $S^2 \times S^1$ when $-w+u = 0$,
because $L_2$ is the trivial knot.

 Since a circle $s_2$ of the surgery slope $-w+u$ bounds
a meridian disc of the Dehn filled solid torus,
a band sum $L'_1$ of $L_1$ and $s_2$ is isotopic to $L_1$
after the Dehn surgery in $L([2w+1, 2u+1])[\infty,-w+u]$.
 We consider the band as illustrated in Figure \ref{fig:BandSum2}
(, where $w = 2, u = -4$ and $s_2$ is of slope $-w+u = -2-4=-6$).
 For this band,
$L'_1$ can be embedded into the torus $\partial N(L_2)$
so that $L'_1 = -2 L -(2u+1) M$.
 It is a $(2, 2u+1)$-cable of the trivial knot $L_2$ in $S^3$,
and a $(2u+1, 2)$-cable of the trivial knot $C'$.
 Note that $L'_1 \cdot M = -2$ and $L'_1 \cdot L = 2u+1$.

 We can obtain the knot $L'_1$ in $L([2w+1,2u+1])[\infty,-w+u]$
by $(-w+u)$-Dehn surgery on $L_2$.
 The new meridian slope is $L+(-w+u)M$.
 Since 
$L'_1 \cdot (L+(-w+u)M) = (-2L -(2u+1)M) \cdot (L+(-w+u)M) = -2(-w+u)+(2u+1) 
= 2w+1$
and 
$L'_1 \cdot M = (-2L -(2u+1)M) \cdot M = -2$,
and since
$M \cdot (L+(-w+u)M) = -1$,
$L'_1$ is a $(2w+1, -2)$-cable of the core $C$ of the filled solid torus.

 We consider the case of $-w+u = \pm 1$,
where $L([2w+1,2u+1])[\infty,-w+u]$ is the $3$-sphere.
 Set $\epsilon = L \cdot (L + (-w+u)M) = -w+u = \pm 1$,
which is the algebraic intersection number
of the meridian of $N(C')$ and that of $N(C)$.
 Hence $L'_1$ is the 
$(L'_1 \cdot (L+(-w+u)M), - \epsilon L'_1 \cdot L) 
= (2w+1, \mp (2u+1))$-torus knot.
\end{proof}

\begin{figure}[htbp]
\begin{center}
\includegraphics[width=8cm]{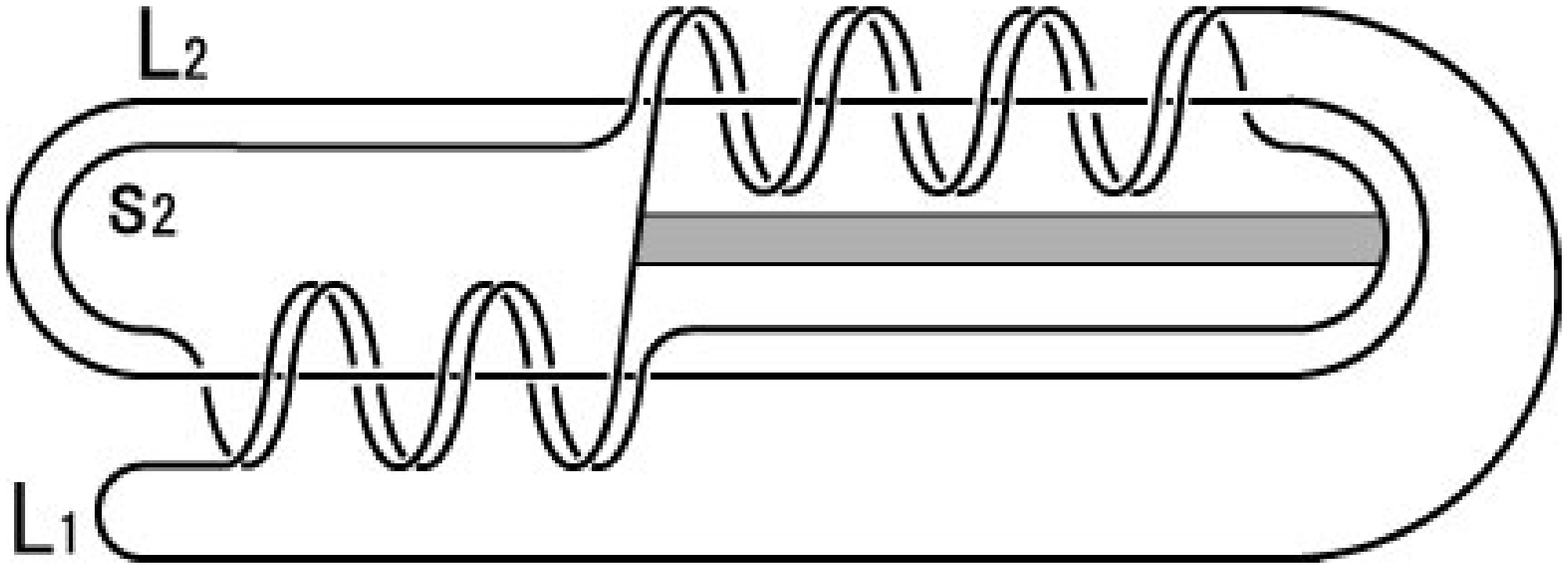}
\end{center}
\caption{}
\label{fig:BandSum6}
\end{figure}

\begin{figure}[htbp]
\begin{center}
\includegraphics[width=5cm]{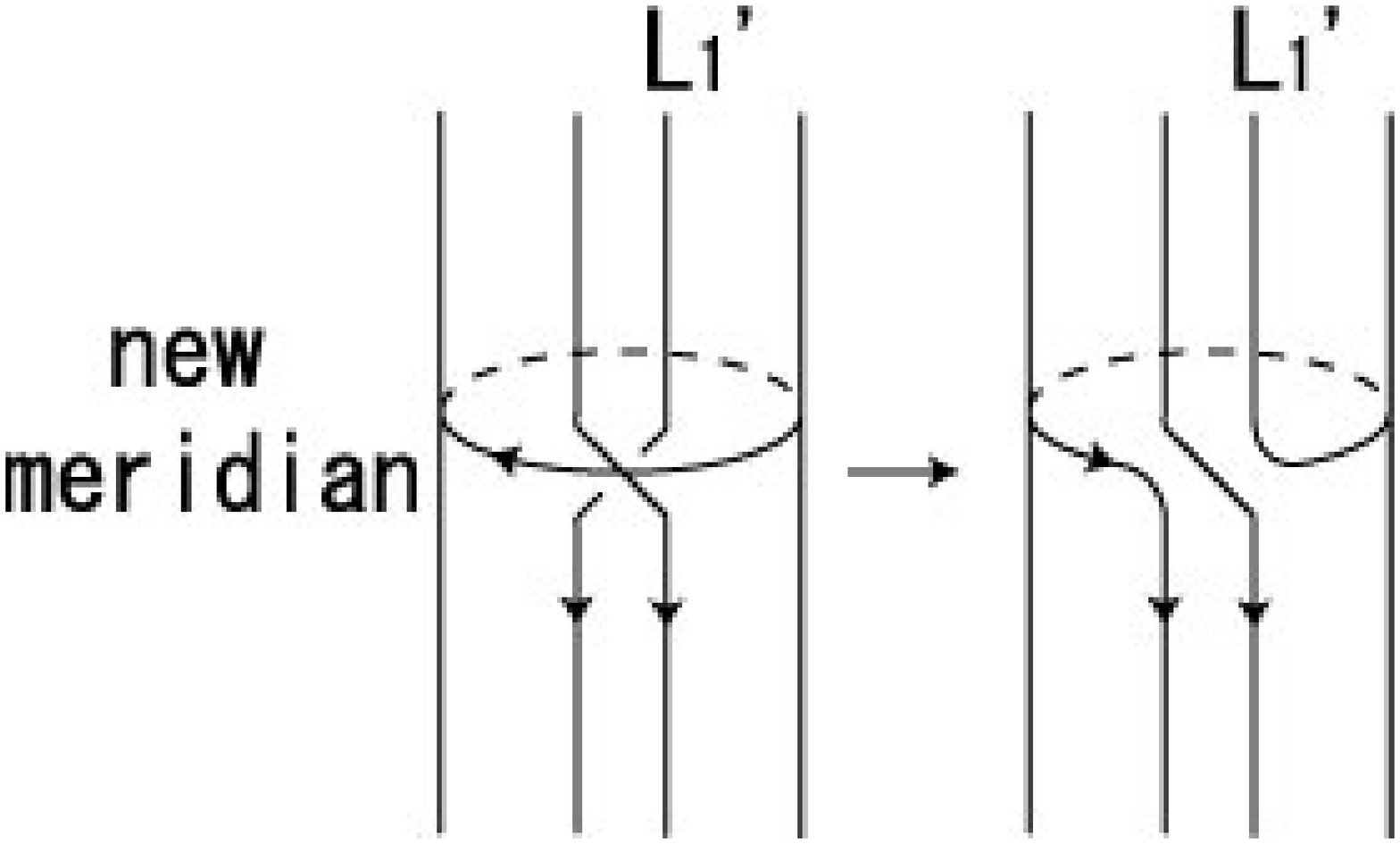}
\end{center}
\caption{}
\label{fig:move}
\end{figure}

\begin{proof}
 We prove Theorem \ref{thm:TorusCable} (2), \ref{thm:torus} (2)
and \ref{thm:cable}.
 Similar argument as that for Theorem \ref{thm:TorusCable} (1)
shows the former half of Theorem \ref{thm:TorusCable} (2).

 The surgery slope is $-w+u+\epsilon$ with $\epsilon = \pm 1$.
 The band sum as in Figure \ref{fig:BandSum6}
(,where $w = 2, u = -4$ and the surgery slope is $-w+u-1 = -2-4-1 = -7$,)
gives the knot $L'_1$.
 For some integer $m$, the knot $L'_1 \subset S^3$
is a $(2, m)$-cable of a knot $K$ in $\partial N(L_2)$,
where $K = -L -\frac{(2u+1) + \epsilon}{2} M$.
 Hence $K$ is a $(1, \{(2u+1) + \epsilon \} /2)$-cable of $L_2$,
and a $(\{(2u+1) + \epsilon \} /2, 1)$-cable of $C'$.
 Note that $K \cdot L = \{ (2u+1) + \epsilon \} /2$,
and $K \cdot M = -1$.

 As in Figure \ref{fig:BandSum6},
$L'_1$ is the $(2,2u+1)$-torus knot in $S^3$ before surgery
if we ignore $L_2$,
and hence $L'_1$ can be isotoped into the torus $\partial N(L_2)$
except near a single $-\epsilon$ crossing
because $(2u+1) - 2 \frac{(2u+1) +\epsilon}{2} = -\epsilon$.

 After the Dehn surgery,
$L + (-w+u+\epsilon) M$ is the new meridian slope of the filled solid torus.
 Hence $K$ winds
$K \cdot (L + (-w+u+\epsilon) M)
= (-L - \frac{(2u+1) +\epsilon}{2} M) \cdot (L + (-w+u+\epsilon) M)
=  -(-w+u + \epsilon) + \{(2u+1) + \epsilon \} /2
= \{(2w+1) - \epsilon \}/2$ times longitudinally
in the filled solid torus.
 Since $K \cdot M = -1$ and $M \cdot (L + (-w+u + \epsilon)M) = -1$,
$K$ is a $(1, -\{ (2w+1) -\epsilon \} /2 )$-cable of $C$.

 $|\{ (2u+1) + \epsilon \} /2|$ is never $0$.
 It is equal to $1$,
if and only if $(u,\epsilon) = (1,-1)$ or $(-2,1)$.
 $|\{ (2w+1) - \epsilon \} /2|$ is never $0$.
 It is equal to $1$
if and only if $(w, \epsilon) = (1,1)$.
 Hence the companion knot $K$ is a core knot,
when and only when $(w, \epsilon)=(1,1)$ or $(u, \epsilon)=(1,-1), (-2,1)$.
 (We consider these cases later.)


 Otherwise, $K$ is a non-trivial non-core torus knot,
and $L'_1$ is a cable knot.
 When $-w+u + \epsilon = \sigma$ with $\sigma = \pm 1$,
the surgery yields the $3$-sphere.
 The case of $(\epsilon, \sigma) = (-1,-1), (-1,1)$
are the mirror images of $(\epsilon, \sigma) = (1,1), (1,-1)$ respectively
with $w$ and $u$ exchanged.
 Hence we consider the case of $\epsilon = 1$ only.
 At this time, the new meridian of the filled solid torus is $L + \sigma M$,
and $K = -L -(u+1)M$.
 The algebraic intersection number with the meridian $L$ of $C'$ is
$L \cdot (L + \sigma M) = \sigma$.
 Since
$K \cdot L 
= (-L - (u+1) M) \cdot L
= u+1$
and 
$K \cdot (L + \sigma M)
= (-L - (u+1) M) \cdot (L + \sigma M)
= (u+1)-\sigma$,
$K$ is the $(\sigma(u+1-\sigma), -(u+1)) =(\sigma(u+1)-1, -(u+1))$-torus knot
in $S^3$.
 The slope of the cabling annulus of $K$ is
$-(\sigma(u+1)-1)(u+1) = -\sigma(u+1)^2 + (u+1)$.
 Since $L'_1$ has a single $-$ crossing in the diagram on $\partial N(L_2)$,
$L'_1$ is the
$(2, 2\{ -\sigma(u+1)^2 + (u+1) \} -1)= (2, -2\sigma(u+1)^2 + 2u+1)$-cable
of $K$.
 Thus Theorem \ref{thm:cable} follows.

\vspace{3mm}

 We consider the case $(w,\epsilon) = (1,1)$.
 (The case $(u,\epsilon) = (1,-1)$ is the mirror image of this case.)
 The surgery slope is $-w+u+\epsilon=u$.
 Recall that $L'_1$ is on the boundary torus of the filled solid torus
except near the $-$ crossing.
 We can move the knot $L'_1$ entirely into the boundary torus
by isotoping a very short under path near the crossing
along a meridian disc of the filled solid torus.
 See Figure \ref{fig:move}.
 Note that the algebraic intersection number
between $K$ and the new meridian is
$(-L - \frac{(2u+1) + 1}{2} M) \cdot (L + u M)
= - u + \frac{(2u+1) + 1}{2}
= 1$ rather than $-1$.
 Hence $L'_1$ is in the position of
$2K - (L + uM)
= 2(- L - \frac{(2u+1) + 1}{2} M) - (L + uM)
= -3L - (3u+2)M$.
 It is a $(3u+2, 3)$-cable of $C'$.
 Since
$(-3L -(3u+2)M ) \cdot (L + uM)
= -3u + (3u+2) = 2$,
$(-3L -(3u+2)M) \cdot M = -3$,
and since $M \cdot (L+uM) = -1$,
$L'_1$ is a $(3,-2)$-cable of $C$.
 When $u = 1$, $L'_1$ is a torus knot in $S^3$.
 At this time, the new meridian of the filled solid torus is $L+M$,
$L'_1 \cdot L = 5$, $L'_1 \cdot (L+M) = 2$.
 Since $L \cdot (L+M) = +1$,
$L'_1$ is the $(2,-5)$-torus knot.
 Thus we obtain Theorem \ref{thm:torus} (2).

 We consider the case $(u,\epsilon) = (-2,1)$.
 At this time, the surgery slope is $-w+u+1=-w-1 (\le -2)$.
 (Such a surgery never yields $S^3$.)
 Hence $K= -L -\frac{(2(-2)+1)+1}{2}M = -L +M$.
 $L'_1$ is on the boundary torus of the filled solid torus
except near the $-$ crossing.
 We can move the knot $L'_1$ entirely into the boundary torus
by isotoping a very short over path near the crossing
along a meridian disc of the complementary solid torus.
 Note that the algebraic intersection number
between $K$ and the meridian $L$ of $C'$ is
$(-L + M) \cdot L = -1$ rather than $+1$.
 Hence $L'_1$ is in the position of
$2K - L = 2(- L +M) - L = -3L + 2M$.
 It is a $(2, -3)$-cable of $C'$.
 Since
$(-3L +2M ) \cdot (L + (-w-1)M) = -3(-w-1) -2 = 3w+1$,
$(-3L +2M) \cdot M = -3$,
and since $M \cdot (L+(-w-1)M) = -1$,
$L'_1$ is a $(3w+1, -3)$-cable of $C$.
\end{proof}

\section{Proof of Theorem\ref{thm:satellite}}


 The last sentence of Theorem \ref{thm:satellite} 
follows from Theorem \ref{thm:TorusCable}.
 We prove the second sentence in this section.

\begin{lemma}\label{lem:EmptyOnL2}
 Let $F$ be a connected surface carried by a branched surface $\Sigma_{\lambda}$
corresponding to a minimal edge-path $\lambda$ for $L(p/q)$.
 Suppose that $F$ has boundary circles on $\partial N(L_1)$
and does no circles on $\partial N(L_2)$.
 Then $\lambda$ is composed of $2m$ edges of label $B$.
 The weights of $F$ are $\alpha = 2$, $\beta = 0$,
the slope of the boundary circles is integral,
and the number of boundary circles is two.
 Moreover, $\chi(F) = -2m+2$ and genus\,$(F) = m-1$.
\end{lemma}

\begin{remark}
  In \cite{S},
Toshio Saito obtained a similar result to this lemma
without using the result of Floyd and Hatcher.
\end{remark}

\begin{proof}
 Since $F$ has not boundary circles on $\partial N(L_2)$,
it cannot be partially carried by $\Sigma_A, \Sigma_C, \Sigma_D$,
but only by $\Sigma_B$ with $\beta = 0$ and $t = \infty$.
 Hence $\Sigma_{\lambda}$ must be composed only of copies of $\Sigma_B$,
and the corresponding minimal edge-path only of edges with label $B$.

 Because the bottom arcs of $\Sigma_B$ are loops,
and since the top arc of $\Sigma_B$ is a single arc,
two consecutive copies of $\Sigma_B$ in $\Sigma_{\lambda}$
are glued along top arcs of both copies, or along bottom ones.
 $\Sigma_B$ has a single branch locus
which is the bottom arc
of the square sector (the line of slope $\infty$) $\times [1/2,1]$.
 We can split $\Sigma_{\lambda}$
along the copies of the square sectors as above,
so that it still carries $F$ after the splittings.
 In fact, such splittings deform $\Sigma_{\lambda}$
into a connected orientable surface $\Sigma$ with no branch locus.
 Hence a surface carried by $\Sigma_{\lambda}$
is the union of $\alpha/2$ parallel copies of $\Sigma$.
 Because $F$ is connected,
it is isotopic to $\Sigma$, and $\alpha = 2$.
 Since $\Sigma$ forms a closed orientable surface when glued along $L_1$,
the slope of the boundary circles is integral,
and the number of boundary circles is two.

 The top and the bottom of $\Sigma_{\gamma}$
are copies of the top of $\Sigma_{\gamma}$
because the bottom of $\Sigma_{\gamma}$ is a union of loops.
 Hence $\Sigma_{\lambda}$ consists of $2m$ copies of $\Sigma_B$
for some positive integer $m$.
 Then $\chi(F) = 2m \frac{\alpha}{2} -(2m-1) \alpha = -2m+2$,
and the genus is $\{ -(-2m+2)-2 + 2 \}/2 = m-1$.
\end{proof}

\begin{proof}
 We prove the second sentence in Theorem \ref{thm:satellite}.
 If $L(p/q)[\gamma]$ is a prime satellite knot,
then its exterior contains an essential torus.
 We assume that the Dehn surgery is performed on $L_1$ rather than $L_2$.
 Then the exterior of the link contains
an essential punctured torus $T$
which has boundary circles only on $\partial N(L_1)$.

 $T$ is carried by a branched surface $\Sigma_{\lambda}$
as in the conclusion of Lemma \ref{lem:EmptyOnL2}
with $m=2$.
 Since $L(p/q) \cong L(p'/q')$ if $p/q - p'/q' \in {\Bbb Z}$,
we can assume that the first edge is from $1/0$ to $0/1$.
 Thus $L(p/q) \cong L([2w, u, 2v])$ for some non-zero integers $w, u, v$.
 If $|w| = 1$ or $|v| = 1$, then $\lambda$ is not minimal,
which is a contradiction.
 If $|u| = 1$, then $L(p/q) \cong L([r,s])$ for some odd integers $r, s$.

 For $L([r,s])$,
minimal edge-paths composed only of edges of label $B$
are $c_5, c_7, d_5$ or $d_7$.
 Substituting $0$ for $\beta$
in the expressions of the boundary slopes for the surfaces
corresponding $c_{25}$, $d_{25}$, $d_{27}$,
we obtain the desired slopes $(-w+u \pm 1, 1)$.
\end{proof}

\begin{figure}[htbp]
\begin{center}
\includegraphics[width=6cm]{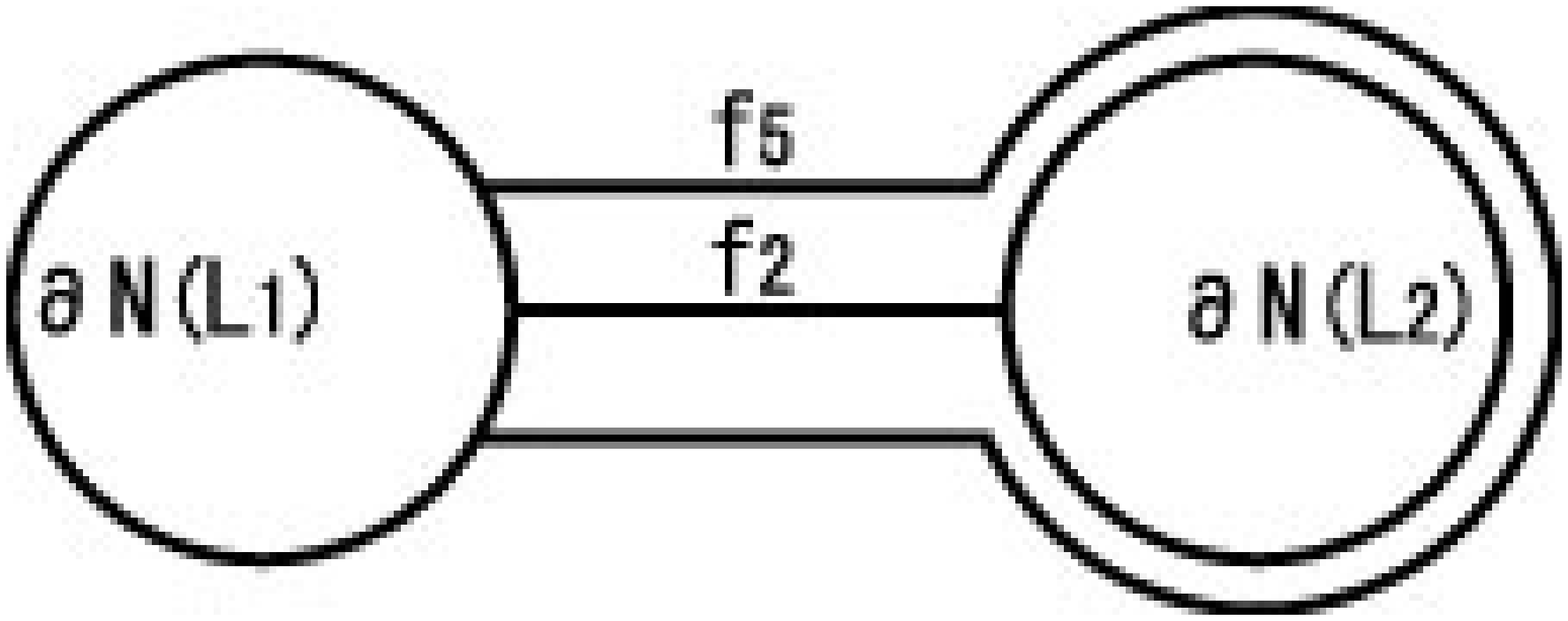}
\end{center}
\caption{}
\label{fig:Double}
\end{figure}

\begin{remark}
 The surface $f_5$ (resp. $f_7$)
corresponding to $c_5$ or $d_5$ (resp. $c_7$ or $d_7$)
with $\alpha = 2$
are obtained
from the double of the surface $f_2$ corresponding to $c_2$ or $d_2$
with $\beta = 1$, $n = 1$ (resp. $n=0$)
by tubing operation on $\partial N(L_2)$
as shown in Figure \ref{fig:Double}.
 Note that the double of $f_2$ gives the cabling annulus
when its boundary circles on $\partial N(L_1)$
are capped off with discs after the Dehn surgery.
\end{remark}

\setcounter{section}{0}
\renewcommand{\thesection}{\Alph{section}}

\section{Calculation for $d_{26}$}

 When $g'=0$, we have the following three equations,
where we set $m = w-u-1 = w+u' (\ge 2)$.

\noindent
$b_1 + b_2 = (m-1)(\alpha - \beta) + \beta + 2 \cdots$ (i)
\newline
$b_1 = GCM(\alpha, m\beta) = GCM(\alpha, m(\alpha-\beta)) \cdots$ (ii)
\newline
$b_2 = GCM(\beta, m\alpha) = GCM(\beta, m(\alpha-\beta))\cdots$ (iii)

 Since $1 < t \le \infty$, it holds $\alpha > \beta \ge 0$.
 When $\beta = 0$, the surface has boundary circles of the meridional slope,
which is not derived from non-trivial Dehn surgery.
 Hence we can assume $\alpha > \beta > 0$.
 Thus $b_1 = GCM(m\beta, \alpha)$ (resp. $b_2 = GCM(m\alpha, \beta)$)
is exact divisor of $\alpha$ (resp. $\beta$)
smaller than or equal to $\alpha$ (resp. $\beta$).

 When $b_1 \le \alpha/2$ and $b_2 \le \beta/2$,
we have $(\alpha + \beta)/2 \ge (m-1) (\alpha -\beta) + \beta + 2$ from (i).
 This implies $(\alpha + \beta)/2 \ge 1 \cdot (\alpha -\beta) + \beta + 2$
since $m \ge 2$.
 Then it follows $\beta -4 \ge \alpha$,
which contradicts $\alpha > \beta$.
 Hence we can assume that $b_1 = \alpha$ or $b_2 = \beta$.

\vspace{3mm}

{\bf The case of $b_1 = \alpha$}.
\newline
 Since $\alpha = b_1 = GCM(\alpha, m(\alpha-\beta))$,
we can see that $\alpha$ is an exact divisor of $m(\alpha-\beta)$,
and there is a positive integer $k$
with $m(\alpha-\beta)= k \alpha \cdots$ (ii)$'$.

 On the other hand, from (iii),
there is a positive integer $h$
with $h b_2 = m(\alpha -\beta) \cdots$ (iii)$'$.
 There are $2$ cases, $h=1$ and $h \ge 2$.

 When $h=1$, we have $b_2 = m(\alpha-\beta) \cdots$ (iii)$''$.
 In the equation (i),
we substitute $\alpha$ for $b_1$ and $m(\alpha-\beta)$ for $b_2$,
to obtain $\alpha - \beta = 1$.
 Then $b_2 = m$ from (iii)$''$, $m = k\alpha$ from (ii)$'$.
 So, again from (iii), $m = b_2 = GCM(\beta, m(\alpha-\beta))$,
which implies that $m$ is an exact divisor of $\beta$ and $m \le \beta$.
 Thus $\beta \ge m = k\alpha \ge \alpha$,
contradicting to $\alpha > \beta$.
 This implies that $h=1$ is inadequate.

 Hence we can assume $h \ge 2$.
 Substituting $\alpha$ for $b_1$ and $\frac{m(\alpha-\beta)}{h}$ for $b_2$
in the equation (i),
we obtain $( m \frac{h-1}{h} -2)(\alpha-\beta) = -2 \cdots$ (i)$'$.
 Comparing signs of both sides,
we can see $m \frac{h-1}{h} -2 < 0$,
and a short calculation gives $m < \frac{2h}{h-1} (\le 4)$.
 Since $h \ge 2$ and $m \ge 2$,
either $m = 2$, or $h=2$ and $m=3$.

 We first consider the case of $m=2$.
 Substituting $2$ for $m$ in (i)$'$, we obtain 
$\alpha -\beta =  h \cdots$ (i)$''$.
 Hence, from (iii)$'$, $b_2 = m = 2 \cdots$ (iii)$'''$.
 From (ii), $\alpha = b_1 = GCM(\alpha, 2\beta)$, which implies
that there is a positive integer $v$ with $v \alpha = 2 \beta$.
 Substituting $\beta + h$ for $\alpha$ by (i)$''$,
we obtain $(2-v)\beta = hv$.
 Hence $2-v > 0$, which implies $v=1$ and $\beta =h$.
 Then $\alpha = 2h$ from (i)$''$.
 From (iii) and (iii)$'''$,
$2 = b_2 = GCM(\beta, m\alpha) = GCM(h, 2 \cdot 2h) = h$.
 Thus we have $\alpha = 4, \ \beta=2, \ w + u' = m = 2$,
which is the unique solution of the equations for $d_{26}$.

 We consider the case where $h=2$ and $m=3$.
 In this case, we obtain from (i)$'$ $\alpha-\beta=4 \cdots$ (i)$'''$.
 Hence, from (iii)$'$, $b_2 = 6$.
 Then, from (iii), $6 = b_2 = GCM(\beta, 12)$,
which implies that $\beta$ is a multiple of $6$.
 From (ii), $\alpha = b_1 = GCM(\alpha, 12)$.
 This implies $\alpha$ is an exact divisor of $12$.
 Hence $\alpha = 12$ and $\beta = 6$ because $12 \ge \alpha > \beta \ge 6$.
 This contradicts (i)$'''$.

\vspace{3mm}

{\bf The case of $b_2 = \beta$}.
\newline
 The equation (ii) implies that
there is a positive integer $h$
with $h b_1 = m(\alpha-\beta) \cdots$ (ii)${}^*$.

 We first consider the case of $h \ge 2$.
 Substituting $\beta$ for $b_2$
and $\frac{m(\alpha-\beta)}{h}$ for $b_1$ in the equation (i),
we obtain $(m \frac{h-1}{h} - 1)(\alpha-\beta) = -2$.
 Comparing signs of both sides,
we can see $m \frac{h-1}{h} -1 < 0$.
 Hence $m < \frac{h}{h-1} \le 2$,
which contradicts $m \ge 2$.

 Hence we can assume $h = 1$.
 Then $b_1 = m(\alpha-\beta)$ from (ii)${}^*$.
 Substituting $\beta$ for $b_2$
and $m(\alpha-\beta)$ for $b_1$ in (i),
we obtain $\alpha-\beta = 2 \cdots$ (i)${}^*$.
 Then $b_1 = 2m$ again from (i).
 From (ii), $2m = b_1 = GCM(\alpha, m(\alpha-\beta))$,
and hence there is a positive integer $k$
with $\alpha = 2mk \cdots$ (ii)${}^{**}$.
 Substituting $\beta+2$ for $\alpha$ by (i)${}^*$,
we have $\beta + 2 = 2mk \cdots$ (ii)${}^{***}$.
 Then from (iii), $\beta = b_2 = GCM(2mk-2, 2m) = GCM(2, 2m)$,
which implies either $\beta = 1$ or $\beta = 2$.
 $\beta = 1$ contradicts (ii)${}^{***}$,
and hence $\beta = 2$.
 Thus $\alpha = 4$ from (i)${}^*$,
and $m = 2$ from (ii)${}^{**}$ and $m \ge 2$.
 This is the unique solution for the equations for $d_{26}$.

\bibliographystyle{amsplain}


\medskip

Hiroshi Goda:
Department of Mathematics,
Tokyo University of Agriculture and Technology,
Koganei, Tokyo, 184-8588, Japan.
goda@cc.tuat.ac.jp

Chuichiro Hayashi:
Department of Mathematical and Physical Sciences,
Faculty of Science, Japan Women's University,
2-8-1 Mejirodai, Bunkyo-ku, Tokyo, 112-8681, Japan.
hayashic@fc.jwu.ac.jp

Hyun-Jong Song:
Division of Mathematical Sciences, 
Pukyong National University, 
599-1 Daeyondong, Namgu, Pusan 608-737, Korea. 
hjsong@pknu.ac.kr

\end{document}